\documentclass[10pt]{article}

\usepackage{amsmath,amsthm,amssymb, epsf,
 amscd,float, setspace,latexsym,amsmath,graphics, 
amsfonts,url,xcolor}
\usepackage{physics}
\usepackage{derivative}
\usepackage{authblk}
 
\usepackage{isomath,amsbsy,euscript,alltt,stmaryrd,enumerate,wrapfig}
\usepackage{mathrsfs} 
\usepackage{ dsfont }
\usepackage[pdftex]{graphicx}
\usepackage{booktabs}
\usepackage{float}
\usepackage{texdraw}
\usepackage{a4wide}
\usepackage{epsfig}
\usepackage{caption}
\usepackage{subcaption}

\usepackage{geometry}
 \def\bn{\mathbf n}

\newcommand{\partderiv}[2]{\ensuremath{\frac{\partial #1}{\partial #2}}}
\newcommand{\partderivtwo}[2]{\ensuremath{\frac{\partial^2 #1}{\partial #2^2}}}

\newcommand{\partderivmix}[3]{\ensuremath{\frac{\partial^2 #1}{\partial #2
\partial #3}}}

\newcommand{\vph}{\vphantom{A^{A}_{A}}}

\newtheorem{corollary}{Corollary}

\newcommand{\bp}{\mathbf p }
\newcommand{\bq}{\mathbf  q}

\newtheorem{theorem}{Theorem}

\newtheorem{proposition}{Proposition}
\newtheorem{remark}{Remark}

\newcommand{\bx}{{\mathbf x}}

\newcommand{\Ccchamber}{C_{\textrm{\begin{tiny}C\end{tiny}}}^{\textrm{\begin{
tiny}II\end{tiny}}}}

\newcommand{\Chcominus}{C_{{\textrm{\begin{tiny}H\end{tiny}}}{\textrm{\begin{
tiny}CO\end{tiny}}}_3^{\textrm{\begin{tiny}-\end{tiny}}}}}

\newcommand{\Jhcominus}{J_{{\textrm{\begin{tiny}H\end{tiny}}}_2{\textrm{\begin{
tiny}CO\end{tiny}}}_3^{-}}^{\textrm{\begin{tiny}II\end{tiny}}}}

\newcommand{\bg}{\mathbf g}

\usepackage{amsmath}
\usepackage{amsbsy}
\usepackage{amssymb}
\usepackage{amscd}
\usepackage{amsfonts}
\usepackage{comment}

\newcommand{\gbext}{{\mathbf g}_{\tiny{\text{ext}}}}
\newcommand{\gbelast}{{\mathbf g}_{\tiny{\text{elast}}}}
\newcommand{\Gelast}{G_{\tiny{\text{elast}}}}
\newcommand{\Gext}{G_{\tiny{\text{ext}}}}

\newcommand{\beq}{\begin{equation}}
\newcommand{\eeq}{\end{equation}}
\newcommand{\beqs}{\begin{eqnarray}}
\newcommand{\eeqs}{\end{eqnarray}}

\newcommand{\psiL}{\psi_{\tiny  \textrm{LdG}}}

\newcommand{\psiB}{\psi_{\tiny \text{Bulk}}}
\newcommand{\psiE}{\psi_{\tiny  \text{E}}}

\numberwithin{equation}{subsection}

\title{Modeling and Analysis of Cooling-Induced Isotropic–Nematic Transitions in Pharmaceutical Liquid Crystals}
\author[1]{Pei Liu}
\author[2]{Debanjan Banerjee}
\author[3]{Daniel Siebel-Cortopassi}
\author[2,4]{Sichen Song}
\author[2]{M. Carme Calderer}

\affil[1]{Department of Mathematics, The Pennsylvania State University\\
224 McAllister Building, University Park, PA 16802, USA\\
\texttt{pul21@psu.edu}} 
\affil[2]{{School of Mathematics, University of Minnesota\\
206 Church Street SE, Minneapolis, MN 55455, USA}\\
\texttt{dbanerje@umn.edu}, \texttt{song0357@umn.edu}, \texttt{calde014@umn.edu}}
\affil[3]{{Department of Mathematics and Systems Engineering, Florida Institute of Technology\\
150 W. University Boulevard, Melbourne, FL 32901, USA}\\
\texttt{dsiebelcorto2023@my.fit.edu}}
\affil[4]{School of Pharmacy, University of Wisconsin--Madison\\
777 Highland Avenue, Madison, WI 53705, USA\\
\texttt{sichen.song@wisc.edu}}

\begin{document}
\maketitle
\graphicspath{{./Fig/}{./figures/}{./Figs/}}

\begin{abstract}

This paper studies the phase transition from an isotropic to a nematic liquid crystal driven by controlled cooling, motivated by the manufacturing of next-generation pharmaceutical drugs. Using the non-isothermal Ericksen model for liquid crystals with a variable degree of orientation, we formulate the problem as a stiff, quasilinear system of parabolic partial differential equations. We establish a consistent thermodynamic framework, prove the well-posedness of the system, and conduct a rigorous asymptotic and numerical scaling analysis of the moving phase boundary in a cylindrical cell subjected to bottom-wall convective cooling. Our analysis identifies a critical sample thickness that separates two distinct physical regimes. Below this threshold, the transition follows the cooling dynamics and the material remains quasi-isothermal. Above the thickness threshold, bulk thermal diffusion  governs the dynamics, imposing a dominant thermal delay characterized by a concave  quadratic  propagation front. These results demonstrate that scaling from laboratory characterization to industrial production fundamentally alters the spatial thermal history of the material, providing a predictive mathematical framework to optimize cooling protocols for kinetic arrest.
\end{abstract}
\section{Introduction.} \label{sec:Intro}
 
In this article, we study the isotropic--nematic phase transition in a thermotropic liquid crystal subjected to a prescribed cooling process. Our work is motivated by the processing of liquid-crystalline pharmaceutical compounds, for which the temperature history can influence the phases formed and, ultimately, the properties of the resulting material. We formulate a non-isothermal extension of Ericksen's theory of liquid crystals with variable degree of orientation, neglecting macroscopic flow. The resulting evolution couples orientational relaxation and heat conduction. We develop a thermodynamically consistent model, establish well-posedness results for the associated quasilinear parabolic initial-boundary value problem, and examine the influence of sample thickness on cooling dynamics through scaling, asymptotic approximations, and numerical calculations.

The solubility and absorption requirements of orally administered drugs have traditionally motivated the use of crystalline and amorphous solid forms. Liquid-crystalline compounds provide an additional source   to develop pharmaceutical materials; their processing is deeply connected with the phase transition properties of such materials.  In particular, the preparation of a glassy material at room temperature may involve cooling from the high-temperature isotropic melt through nematic and, under suitable conditions, smectic phases. The cooling rate influences which transitions occur before the material becomes kinetically arrested. We focus on itraconazole (ITZ), a thermotropic, glass-forming pharmaceutical compound for which calorimetric and X-ray-scattering measurements provide a useful experimental setting. Upon cooling, ITZ exhibits an isotropic--nematic transition near $\theta_{NI}=363\,\mathrm{K}$ and, under sufficiently slow cooling, a nematic--smectic transition near $347\,\mathrm{K}$. Our present analysis concerns the first of these transitions only. 

Ericksen's model with variable degree of orientation provides the starting point for the governing equations \cite{Ericksen1991}. In the absence of flow, the orientational variables are the unit director field $\bn$, describing a preferred molecular alignment direction, and the scalar order parameter $s$, describing the degree of alignment. The value $s=0$ corresponds to an isotropic distribution, while $s=1$ and $s=-\frac12$ represent limiting states of perfect alignment parallel and perpendicular, respectively, to the director. The physically admissible range is $-\frac12<s<1$. The director is also the central variable of Oseen--Frank theory, whereas the additional scalar field in Ericksen's theory allows the degree of orientational order to vary, as required in the vicinity of an isotropic--nematic transition.

The model is closely related to the Landau--de Gennes theory of the symmetric, traceless order tensor $Q$ under the uniaxial representation
$
Q=s\left(\bn\otimes\bn-\frac13 I\right).
$
We consider both the classical polynomial Landau--de Gennes bulk potential and the singular potential introduced by Ball and Majumdar on the basis of Maier--Saupe mean-field theory. The polynomial potential, although convenient for analysis and computation, is defined for arbitrary order tensors and does not by itself exclude unphysical eigenvalues. On the other hand, the Ball--Majumdar potential incorporates the constraints imposed by the statistical-mechanical definition of $Q$: its entropic contribution diverges as an eigenvalue approaches the boundary of the admissible set. In the uniaxial case, this divergence occurs as $s\to-\frac12$ or $s\to1$, providing an energetic mechanism for upholding the physical bounds on molecular order.

The Ball--Majumdar potential is, however, specified implicitly through a constrained entropy minimization problem, even after reduction to uniaxial symmetry. To make this description more accessible for analysis and computation, we construct an explicit approximation by separating the logarithmic endpoint singularities from a regular remainder and approximating the latter using Chebyshev polynomials. We prove convergence of the potential and its first two derivatives on compact subintervals of the admissible range, together with convergence of the associated coexistence order parameter and interaction coefficient. The approximations retain the endpoint singularities and satisfy the isotropic--nematic coexistence conditions. These properties allow us to preserve the principal physical features of the singular theory while working with an explicit expression for the bulk energy.

The evolution equations for $\bn$ and $s$ follow from the balance of generalized angular momentum after neglecting the small inertial terms associated with second-order time derivatives. They can alternatively be formulated through the principle of minimum dissipation \cite{CaldererGolovatyLinLiu2002}, \cite{CaldererGolovatyLavrentovichWalkington2016}, \cite{TovkachCaldererGolovatyLavrentovichWalkington2016}. To incorporate thermal effects, we introduce the temperature $\theta$ as an additional unknown and specify constitutive relations for the Helmholtz free energy, entropy, internal energy, and heat flux. The balance of energy and the Clausius--Duhem inequality then yield a coupled, thermodynamically consistent evolution system. Changes in orientational order contribute to the heat equation, while the thermal conductivity tensor accounts for different conductivities parallel and perpendicular to the director. The formulation also identifies the dissipative mechanisms associated with orientational relaxation and heat conduction.

Motivated by the experimental arrangement, we consider a sample in a cylindrical cell cooled through its lower base by contact with  a cooling bath, and  with insulated lateral and upper surfaces. Heat exchange with the temperature-controlled reservoir is represented by a Newton--Robin condition, allowing for a finite thermal contact resistance. Consequently, the imposed reservoir temperature need not coincide with the temperature of the liquid crystal at the lower boundary. We prescribe vertical director anchoring on the lateral wall, natural conditions for the director on the two bases, and homogeneous Neumann conditions for the scalar order parameter. These conditions admit a spatially uniform vertical director; for initial data depending only on the axial coordinate, the coupled system admits a corresponding one-dimensional reduction for $s$ and $\theta$. The model is formulated for the laboratory-scale geometry, while thicker samples are considered to assess the implications of changing the thermal length scale. Extrapolation to industrial processing would require additional physical effects beyond those included here.

We establish local existence, uniqueness, and regularity for the coupled initial-boundary value problem, together with separation of the order parameter from the physical bounds and a continuation criterion under the prescribed assumptions. This mathematical framework supports the subsequent study of the reduced equations. Nondimensionalization reveals strongly separated relaxation and heat-transfer time scales. In particular, comparison of the thermal diffusion time across the sample with the characteristic time of external cooling identifies a critical height  
$
d_{\mathrm{crit}}=\left(\frac{\kappa_{\parallel}t_0}{\rho c_{\theta}}\right)^{1/2},
$ 
which is approximately $4.56\,\mathrm{mm}$ for the reference material parameters and cooling time used here. For thinner samples, temperature variations across the cell are relatively small on the bath-cooling time scale; for thicker samples, internal diffusion produces a more pronounced spatial delay. 

Using the polynomial Landau--de Gennes potential for the scaling calculations, we obtain approximate temperature fields in these regimes and examine the evolution of the isotherm $\theta=\theta_{NI}$ at different heights. The resulting arrival-time curves have an approximately concave quadratic dependence on the axial coordinate. This transition-temperature isotherm, although  a useful indicator of the thermal evolution of cooling,  it may not coincide with an actual isotropic--nematic interface: metastability and the kinetics of orientational relaxation may delay the change in $s$ after the coexistence temperature has been crossed. Our calculations therefore distinguish thermal predictions from the full evolution of orientational order.

The paper is organized as follows. Section~B reviews the Landau--de Gennes and Ball--Majumdar descriptions of the isotropic--nematic transition and develops the explicit approximation of the singular potential, including its convergence properties. Section~C presents the non-isothermal Ericksen model, its thermodynamic derivation, and the initial and boundary conditions associated with the experimental geometry. Section~D establishes the well-posedness results for the coupled parabolic system and justifies the symmetry reduction. Section~E develops the nondimensionalization, identifies the relevant thermal time scales, and analyzes the temperature evolution and transition-temperature isotherms through asymptotic and numerical approximations. Section~F discusses the conclusions. Additional arguments and calculations are collected in the Supplementary Material.

\section{Modeling Background.} \label{MathPrelim}
We consider non-isothermal nematic liquid crystals modeled according to the Ericksen theory of liquid crystals with variable degree of orientation \cite{Ericksen1991}, consistent with the  Landau-de Gennes model restricted to uniaxiality. We start with a brief review of the latter that will also  help us set up the notation.  The anisotropy of a nematic liquid crystal is quantified by the traceless symmetric order tensor $\mathbf Q$ \cite{DeGennes1971}: 
\begin{align}
    \mathbf{Q} (\mathbf{x}) :=  \int_{\mathbb{S}^2} (\mathbf{p} \, \otimes \, \mathbf{p} - \frac{1}{3} \mathbf{I}) \, \rho_0 (\mathbf{x}, \mathbf{p}) \, d \mathbf{p}, \label{eq:Qtensor}
\end{align}
the normalized second order moment of the probability density distribution  $\rho_0(\mathbf{x}, \mathbf{p})$ of molecular orientation vectors $\mathbf{p} \in \mathbb{S}^2$ at a point $\mathbf{x}\in\Omega$, domain occupied by the liquid crystal, and  with $\mathbf{p}$ and $\mathbf{-p}$ being equivalent. 
A property of $\mathbf Q$ is the lower bound satisfied by its minimum eigenvalue, $\lambda_{\textrm{min}} (\mathbf{Q} (\mathbf{x})) \geq -\frac{1}{3}.$
\subsection{Static theory}
The total free energy, at temperature $\theta$, is postulated as
\begin{align}\label{eq:total-energy}
 &   I_\theta(\mathbf{Q}, \theta) =\int_\Omega  (\psiB(\mathbf{Q}, \theta) + \psiE (\mathbf{Q}, \nabla \mathbf{Q})\, d \mathbf{x},
\end{align}
where   $\psiB$ and $\psiE$ denote the bulk and elastic contributions, respectively. The elastic contribution, $\psi_E (\mathbf{Q, \nabla Q})$, is  taken as a quadratic, isotropic function of $\nabla \mathbf{Q}$, such as 
\begin{align}
    \psi_E (\mathbf{Q, \nabla Q}) = \sum^4_{i=1} L_i \mathrm{I_i}, \label{eq:phiE}
\end{align}
where $\mathrm{I_i} = \mathrm{I_i} (\mathbf{Q, \nabla Q})$ denote the invariants 
$ \mathrm{I_1} = \mathrm{Q_{ij,k}} \mathrm{Q_{ij,k}}, $ $ \mathrm{I_2} = \mathrm{Q_{ij,j}} \mathrm{Q_{ik,k}}, $ $\mathrm{I_3} = \mathrm{Q_{ik,j}} \mathrm{Q_{ij,k}}, \, $ $ \mathrm{I_4} = \mathrm{Q_{lk}} \mathrm{Q_{ij,l}}, \mathrm{Q_{ij,k}}.$
 In the special case that  $L_2 = L_3 = L_4 = 0$, the latter yields
  $ \psi_E = L |\nabla \mathbf{Q}|^2. $
We will focus on uniaxial liquid crystals for which  the following  order tensor representation holds:
\begin{align}
    \mathbf{Q} := s (\mathbf{n} \otimes \mathbf{n} - \frac{1}{3} \mathbf{I}), \mathbf{n} \in \mathbb{S}^2,\,\,\, \, |\bn|=1, \,\, s\in (-\frac{1}{2}, 1). \label{eq:Quniaxial}
\end{align}
We note that the bulk energy $\psiB(Q, \theta)$ reduces to a function of $s$, in the uniaxial case.   For notational convenience, we will use the same symbol to represent the latter, that is,
$\psiB(s,\theta)$. Plot of thermotropic energy ($W_0$) vs the order parameter ($s$).

\begin{figure}[htbp]
    \centering
 \begin{subfigure}[b]{0.35\textwidth}
        \centering
        \includegraphics[width=\textwidth]{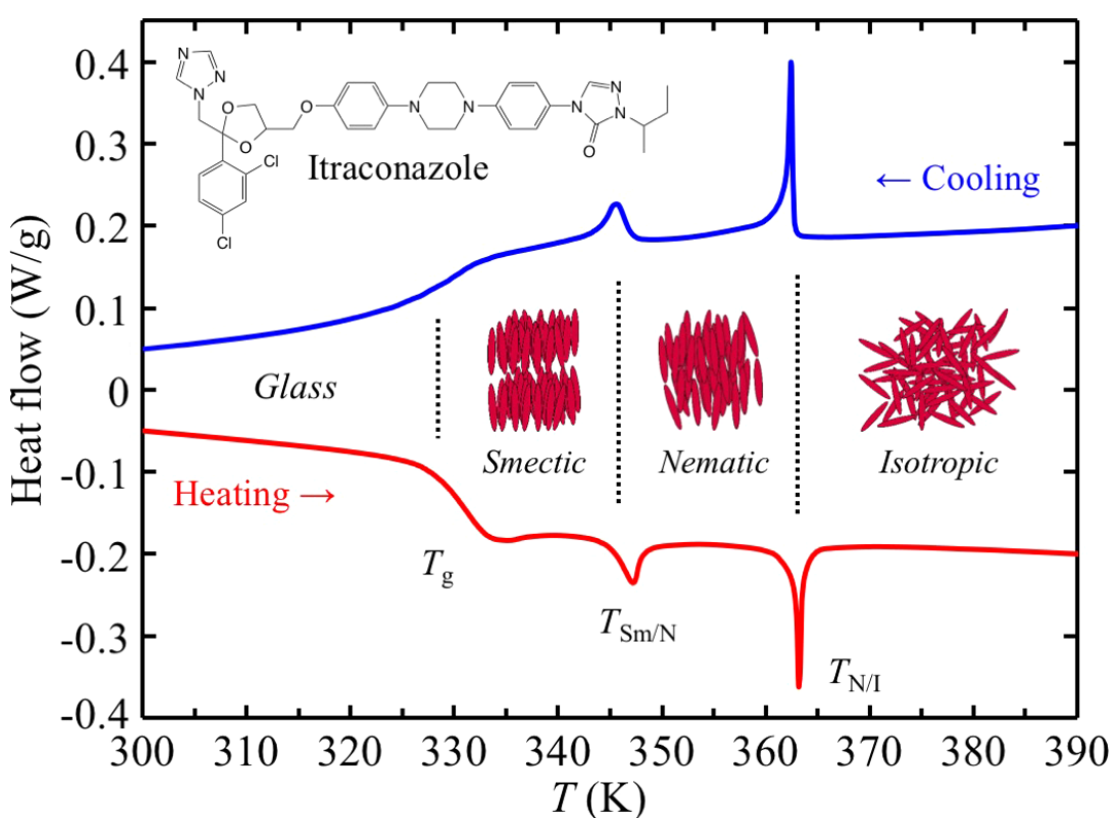}
        \label{fig:heat_flow}
        \caption{\small{Heat flow versus temperature}}
    \end{subfigure}
    \hfill 
    \label{fig:experiment}
    \begin{subfigure}[b]{0.4\textwidth}
        \centering
        \includegraphics[width=\textwidth]{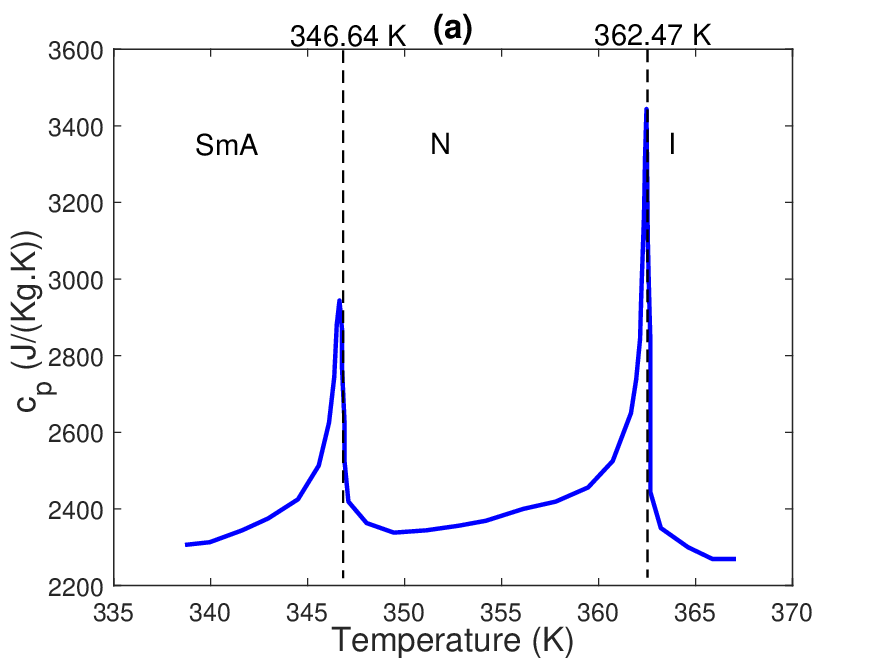}
          \caption{\small{Specific heat versus temperature}} \label{fig:specific_heat}
    \end{subfigure}
    \hfill
    \caption{\textbf{(a)}: Differential Scanning Calorimetry traces of ITZ cooling and reheating at 0.16 K/s. The isotropic-nematic transition ($T_{\textnormal{I/N}}$ = 363 K), the nematic-smectic transition ($T_{\textnormal{Sm/N}}$ = 347 K), and the glass transition ($T_\textnormal{g}$ = 328 K). The upper panel shows the molecular shape of ITZ.
   \cite{Teerakapibal2018}. \textbf{(b)} Variation of $c_p$ with temperature for Smectic A (SmA), Nematic (N), and Isotropic phases (I) \cite{AmponsahEfah2020}. 
    \label{fig:double_well_transition}.}
\end{figure}
\subsection{Landau-de Gennes Model}
In this model, the bulk energy 
   $\psiL$ is  taken as a truncated Taylor expansion about $\mathbf{Q} = \mathbf{0}$ \cite{SchopohlSluckin1987},
\begin{align}
 \psiB(\mathbf{Q}, \theta):=   \psiL (\mathbf{Q}, \theta) = a \operatorname{tr} (\mathbf{Q}^2) + \frac{2b}{3} \operatorname{tr} (\mathbf{Q}^3) + \frac{c}{2} (\operatorname{tr} (\mathbf{Q}^2))^2, \quad  a = \alpha (\theta - \theta^*), \label{eq:psi-LdG}
\end{align}
where $ c > 0, \, \alpha>0$  and $\theta^*$ is the fixed reference temperature at which the isotropic state becomes unstable. 
Existence of minimizer for the Landau de Gennes energy \eqref{eq:total-energy} and \eqref{eq:psi-LdG} subject to Dirichlet boundary conditions, and under appropriate conditions on the coefficients $L_i$ is analyzed in \cite{davis1998finite}. 
It was later established that,  for prescribed $\mathbf{Q}$ uniaxial on the boundary, and for sufficiently small values of $|L_i|$, $1\leq i\leq 3$ (the parameter range of interest), the minimizer $\mathbf{Q}$ is uniaxial, except, perhaps, in a set of measure zero \cite{BallMajumdar2010}.   
We point out that, for $Q$ uniaxial, the Landau-de Gennes energy density \eqref{eq:psi-LdG} reduces to

\begin{equation}\psiL(s, \theta):= \frac{2}{27}(9as^2+ 2bs^3+3cs^4), \quad a = \alpha (\theta - \theta^*).  \label{psi_unaxial}\end{equation}

\noindent {\bf Phase transition from Isotropic to Nematic.\,} 
In order to characterize the transition temperatures and the range of coexistence of both states, we proceed to the stability analysis as shown in \cite{MottramNewton2014}. For this, we calculate the critical points of $\psiL$ and the corresponding values of
the second derivatives $\frac{d^2 \psiL}{d s^2}$,  and subsequently compare the energies of such  points. We arrive at the following \cite{MottramNewton2014}:

\begin{enumerate}
 \item The isotropic state $s=0$ is globally stable for $a > \frac{b^2}{27c}$, metastable for $0 < a < \frac{b^2}{27c}$ and unstable for $a<0$. 
\item  The nematic phase,  $s = \frac{-b + \sqrt{b^2 -24ac}}{4c}$, is globally stable for $a < \frac{b^2}{27c}$ metastable for $\frac{b^2}{27c} < a < \frac{b^2}{24c}$, and undefined  for $a > \frac{b^2}{24c}$.

\item The phase $s = \frac{-b - \sqrt{b^2 -24ac}}{4c}$ is metastable (with $s$ negative) for $a < 0$, unstable for $\frac{b^2}{27c} < a < \frac{b^2}{24c}$, and not defined for $a > \frac{b^2}{24c}$.
\end{enumerate}
Since our goal is to study the phase transition from the isotropic to the nematic upon cooling, the relevant temperature is precisely $\theta_\text{NI}$, the lowest temperature where the isotropic state is still locally stable. In fact, thre
e threshold temperatures, $\theta^*<\theta_\text{NI}<\theta^+$ have been identified from figure \ref{fig:double_well_transition}(a) as follows. (1) $\theta^*$: The low temperature at which the isotropic state becomes unstable, (2) $\theta_\text{NI}$: Nematic-Isotropic transition temperature, that is, the temperature at which the energy of the isotropic and nematic states are exactly equal, and (3) $\theta^+$: The high temperature where the nematic state disappears.  Consequently, there are three threshold values of $a$: (1) $a = 0$, corresponding to $
\theta^*$, (2) $a = \frac{b^2}{27c}$ corresponding to $\theta = \theta_{NI}$, and (3) $a = \frac{b^2}{24c}$, corresponding to $\theta = \theta^+$. If we use the notation  $a = \alpha (\theta - \theta^*) = \alpha \mathrm{\Delta} \theta$, then the threshold temperatures are \cite{MottramNewton2014}: 
\begin{equation}
\theta^+ = \frac{b^2}{24ac} + \theta^*, \quad \theta_{\textnormal{NI}} = \frac{b^2}{27ac} + \theta^*, \quad \theta = \theta^*.\end{equation}
Therefore, for $\psiL$ and depending on the coefficient values $a, b, c$, in an intermediate temperature region ($\theta=\theta_{\textrm{NI}}$), the minimum in the isotropic state and the minimum in the uniaxial nematic states may both be locally stable (see figure \ref{fig:double_well_transition}(b)). In fact, at sufficiently low temperatures ($\theta=\theta^*$), the isotropic state must lose stability, leaving only the uniaxial nematic state stable, that is, only one global minimum close to $s=1$. In contrast, at sufficiently high temperatures ($\theta=\theta^+$), the uniaxial states must lose stability, leaving only the isotropic state stable, that is, only one global minimum at $s=0$. The values of the parameters $b, c$ are given by $b=-0.64\times10^6 \text{N}\text{m}^{-2}, c=0.35\times10^6 \text{N}\text{m}^{-2}$ \cite{MottramNewton2014}. From figure \ref{fig:double_well_transition}(a), the values of the three threshold temperatures are chosen or computed as follows: $\theta_\text{NI}=363\text{K},\, \theta^*=350\text{K},\, \theta^+=\frac{b^2}{24 \alpha c}+\theta^*$, where $\alpha$ is computed from $\theta_\text{NI}=\frac{b^2}{27 \alpha c}+\theta^*$, that is, $\alpha=\frac{b^2}{351c}$.

\subsection{The Maier-Saupe model and the Ball-Majumdar potential}
\begin{figure}[htbp]
    \centering
 \begin{subfigure}[b]{0.40\textwidth}
        \centering
        \includegraphics[width=\textwidth]{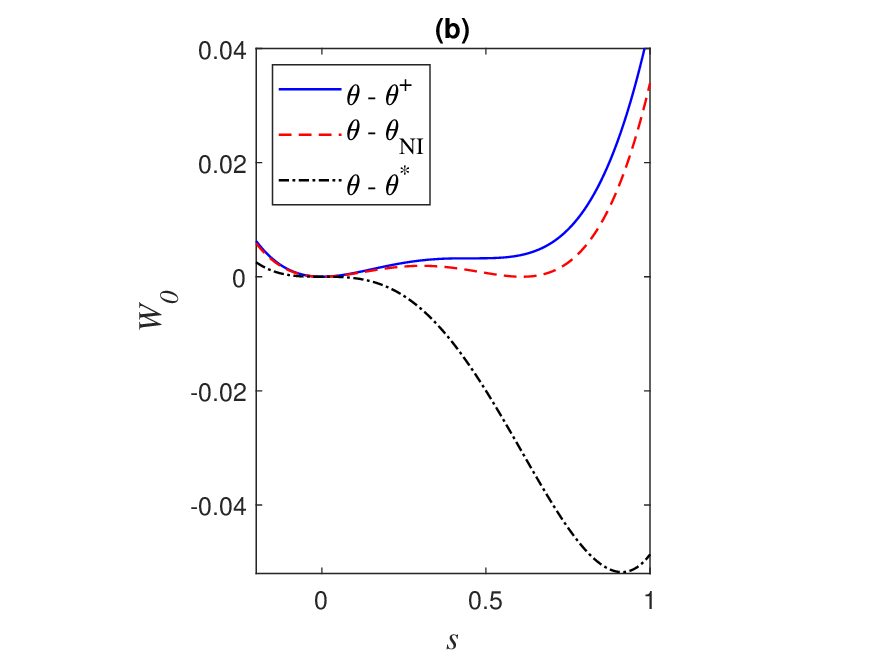}
        \caption{Singular potential}
        \label{fig:double-well}
    \end{subfigure}
    \hfill 
    \begin{subfigure}[b]{0.4\textwidth}
        \centering
        \includegraphics[width=\textwidth]{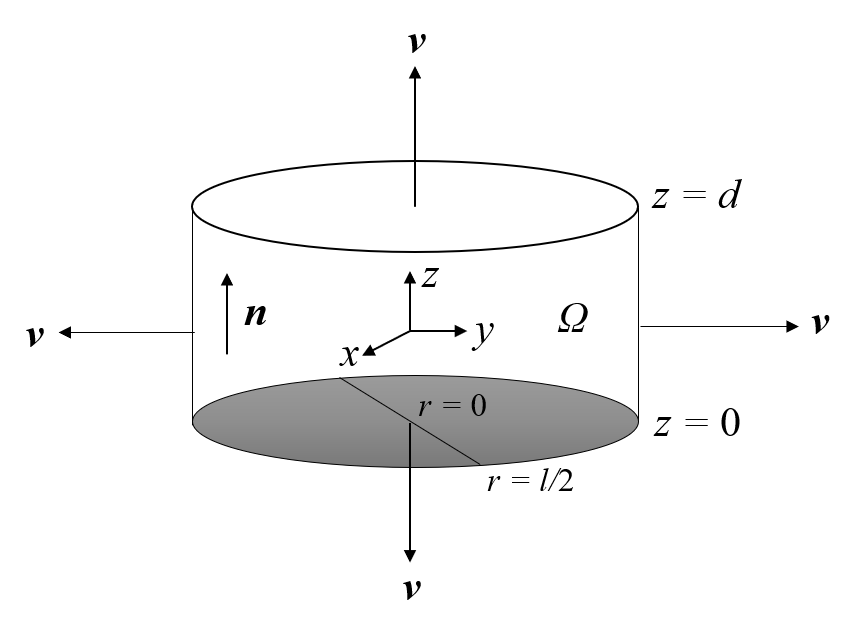}
        \caption{Liquid crystal domain $\Omega$}
        \label{fig:cylinder_domain}
    \end{subfigure}
    \hfill 
    \caption{\small{ \textbf{(a)} Double well potential in the Ball-Majumdar model. \textbf{(b)} Cylindrical domain  of height \text {d} occupied by the liquid crystal. 
    }} 
    \label{fig:combine_profiles}
\end{figure}

Starting with  the Maier-Saupe mean field theory, Ball and Majumdar, following earlier works by Katriel et al., constructed a macroscopic bulk energy for which the physical bounds of the eigenvalues of $Q$ hold \cite{BallMajumdar2010, Katriel1986KKLS,Katriel1983MeanField}. The construction starts with the cross-graining of the Maier-Saupe entropic component of the energy by minimizing it among all probability density function $\rho$ that yield a given order tensor $\mathbf {Q}$, 
\begin{equation} \label{f(Q)}
    f(Q)= \begin{cases}
\inf_{\rho\in\mathcal{} A_Q}\int_{\mathcal S^2}\rho(\bn)\ln \rho(\bn)\,dS \\[1ex]
+\infty, \, \text{otherwise},
    \end{cases}
\end{equation} with
\begin{equation} \label{eqn:AQ}
\mathcal A_ModQ= \left\{\rho: {\mathcal S}^2\longrightarrow \mathbb R,
\rho\geq 0, \quad \int_{{\mathcal S}^2} \rho(\bp)\,d\bp=1, \quad Q=\int_{{\mathcal S}^2}(\bp\otimes\bp-\frac{1}{3} I)\rho(\bp)\, d\bp
\right\}
\end{equation}
With this, the bulk energy in equation \eqref{eq:total-energy} is now taken as 
\begin{equation} 
\psiB(Q, \theta):=    \psi(Q, \theta)= \theta f(Q)- \kappa |Q|^2. \label{psi_MS}
\end{equation}
In  \cite{BallMajumdar2010}, the authors establish existence of minimizer in \eqref{f(Q)} and the asymptotic result $f(Q)\to +\infty$,  as one of the eigenvalues $\lambda_i(Q)\to -\frac{1}{3}$.  Their construction of such a potential, in the uniaxial case,  is highlighted next. 
For $\mathbf{Q}$ as in \eqref{eq:Quniaxial}, the potential in \eqref{f(Q)} admits the form 
\[
f(s)=\frac{1+2s}{3}\mu(s)-\ln \left( Z(\mu(s)) \vph\right),\,\, s\in (-\frac{1}{2}, 1)
\quad \textrm{and}\quad 
Z(\mu)=\int_{S^2}e^{\mu p_3^2}\,dp
=2\pi\int_{-1}^{1}e^{\mu u^2}\,du,
\]
with $\mu(s)$ uniquely determined by the relation
\[
\frac{d}{d\mu}\ln\left( Z(\mu(s)) \vph\right) = \frac{1+2s}{3},
\quad \text{ with} \quad
Z(\mu)=
\begin{cases}
\dfrac{2\pi^{3/2}}{\sqrt{\mu}}\operatorname{erf}(\sqrt{\mu}),
& \mu>0,\\[10pt]
\dfrac{2\pi^{3/2}}{\sqrt{-\mu}}\operatorname{erf}(\sqrt{-\mu}),
& \mu<0 .
\end{cases}
\]
Furthermore, the following asymptotic relations hold, 
\begin{align}
f(s)=-\ln(1-s)+O(1),
\,\,
f(s)=-\tfrac12\ln(1+2s)+O(1), \,\, \label{limits-BM} \\   f(s)= -\ln(4\pi)+\frac52 s^2+\frac{25}{21}s^3+\frac{425}{196}s^4+O(s^5), \nonumber \end{align}
as $s\to 1^-,  s\to-\tfrac12^+ $ and $s\to 0$, respectively. 
With this, the uniaxial free energy in \eqref{psi_MS} becomes 
\begin{equation} 
\psi(s, \theta)= \theta f(s)- \frac{2}{3}\kappa s^2. 
\label{energy_BM_uniaxial}
\end{equation}
We note that, in the uniaxial setting, the Ball--Majumdar bulk potential provides a macroscopic order-parameter representation of the Maier--Saupe mean-field free energy.  In both approaches, the quadratic term represents an intermolecular interaction, with strength $\kappa$, that in view of \eqref{energy_BM_uniaxial},  parametrizes such a family of potentials. It is illustrated in figure \ref{fig:double-well}.

\subsubsection{Characterization of the Bulk Potential in the Isotropic-Nematic Phase Transition}

The Maier-SThe Maier-Saupe potential identifies the approximate value of the order parameter $s_{NI}\approx 0.429$  as associated with the  temperature of transition,  $\theta_{NI}$, between the isotropic and the nematic phases \cite{MaierSaupe1959, MaierSaupe1960}. This property is inherited by the Ball--Majumdar potential whose approach predicts the existence of a transition value of $s$ in the interval $[\frac{1}{4}, \frac{1}{2}]$ \cite{BallMajumdar2010}. Such a value remains the same regardless the temperature of phase transition of a particular liquid crystal. 
aupe potential identifies the approximate value of the order parameter $s_{NI}\approx 0.429$  as associated with the  temperature of transition,  $\theta_{NI}$, between the isotropic and the nematic phases \cite{MaierSaupe1959, MaierSaupe1960}. This property is inherited by the Ball--Majumdar potential whose approach predicts the existence of a transition value of $s$ in the interval $[\frac{1}{4}, \frac{1}{2}]$ \cite{BallMajumdar2010}. Such a value remains the same regardless the temperature of phase transition of a particular liquid crystal. 

A problem of interest in applications is: given the phase transition temperature $\theta_{NI}$ of the liquid crystal under study, determine the corresponding potential $\psi(s, \theta_{NI})$ within the family \eqref{energy_BM_uniaxial}. Its construction  together with the universality of $s_{NI}$ follows from the coexistence conditions, 
\begin{equation}
    \psi(0, \theta_{NI})=\psi(s_{NI}, \theta_{NI}), 
    \qquad \partderiv{\psi}{s}(0, \theta_{NI})=0=\partderiv{\psi}{s}(s_{NI}, \theta_{NI}).
\end{equation}
The first derivative condition follows directly from the construction of $f(s)$; the two remaining  give
\begin{equation}
    \kappa=\frac{3\theta_{NI}}{4 s_{NI}}f'(s_{NI})
    \quad  \text{ and } \quad
    \Phi(s_{NI})=0,
    \qquad \text{ where } \quad
    \Phi(s):= f(s)- f(0)-\frac{1}{2} s f'(s). \label{coexistence-conditions}
\end{equation}
Solving the latter gives the previously indicated value $s_{NI}$ as obtained from the Maier-Saupe and Ball-Majumdar theories. The uniqueness of the positive root follows from standard calculus arguments. 

\subsubsection{Explicit approximation of the Ball--Majumdar potential}
\label{sec:BM-approx}

One disadvantage of the Ball--Majumdar potential is that its dependence
on the scalar order parameter is given only implicitly.  We construct
an explicit approximation which retains both the logarithmic
singularities at the endpoints of the physical interval and the
double-well structure at the isotropic--nematic transition.

We first introduce some notation and properties of polynomial approximation. The following function encodes the singular behavior of the potential \eqref{energy_BM_uniaxial},
\begin{equation}
\label{eq:singular-part}
L(s)
=
-\log(1-s)-\frac12\log(1+2s),
\qquad
-\frac12<s<1,
\end{equation}
so that 
\begin{equation}
\label{eq:regular-remainder}
f(s)=L(s)+R(s),
\qquad
R(s):=f(s)-L(s).
\end{equation}
Fix $\delta>0$ such that  \begin{equation}
s_{NI}\in I_\delta=[a,b]
=
\left[-\frac12+\delta,\,1-\delta\right]. \label{I_delta}
\end{equation}
On $I_\delta$, the function $R$ is regular.
Let us  introduce the affine map $I_{\delta}\longrightarrow [-1, 1]$, 
\begin{equation}
\label{eq:Cheb-map}
x=x(s)
=
\frac{2s-a-b}{b-a},
\qquad
s\in I_\delta,
\end{equation}
 and define
\[
\widehat R_2(x)
=
R''\left(
\frac{a+b}{2}+\frac{b-a}{2}x
\right).
\]
Let
\begin{equation}
\label{eq:Cheb-qm}
q_m(s)
=
\sum_{k=0}^{m}c_kT_k(x(s)),
\end{equation}
where $T_k$ denotes the k-th degree Chebyshev polynomial of the first kind and
\begin{equation}
\label{eq:Cheb-coefficients}
c_k
=
\frac{2-\delta_{k0}}{\pi}
\int_0^\pi
\widehat R_2(\cos\vartheta)
\cos(k\vartheta)\,d\vartheta .
\end{equation}
We denote the Chebyshev approximation error by
\begin{equation}
\label{eq:Cheb-error}
e_m
:=
\|q_m-R''\|_{L^\infty(I_\delta)}.
\end{equation}
The polynomial approximation of the regular part is obtained by
integrating twice:
\begin{equation}
\label{eq:Pm-definition}
P_{m+2}(s)
=
R(0)+R'(0)s
+
\int_0^s(s-\eta)q_m(\eta)\,d\eta,
\end{equation}
and we set
\begin{equation}
\label{eq:fm-definition}
f_m(s):=L(s)+P_{m+2}(s).
\end{equation}
Corresponding to the coexistence function $\Phi$ defined in
\eqref{coexistence-conditions}, set
\begin{equation}
\label{eq:Phi-m}
\Phi_m(s)
=
f_m(s)-f_m(0)-\frac12 s f_m'(s).
\end{equation}
The approximate coexistence value $s_{NI}^{(m)}$ is determined by
\begin{equation}
\label{eq:approx-coexistence}
\Phi_m(s_{NI}^{(m)})=0.
\end{equation}
The corresponding approximating bulk potential is

\begin{equation}
\label{eq:psi-m}
\psi_m(s,\theta)
=
\theta f_m(s)-\frac23\kappa_m s^2, \quad\text{with}\quad  \kappa_m
=
\frac{3\theta_{NI}}{4s_{NI}^{(m)}}
f_m'\bigl(s_{NI}^{(m)}\bigr)
\end{equation}
The following theorem establishes an explicit approximation of the Ball-Majumdar potential and a procedure to calculate it. 
\begin{theorem}
\label{thm:BM-Chebyshev}
 Let $s_{NI}$ be the coexistence order  parameter of the Maier-Saupe theory and  $\psi $ as in \eqref{energy_BM_uniaxial} denote the uniaxial Ball--Majumdar potential, with the interaction parameter $\kappa$ parametrized by $\theta_{NI}$ as in  \eqref{coexistence-conditions}. Let $\delta>0$ and the interval $I_{\delta}$ be as in \eqref{I_delta}. Let the singular and regular parts of the potential, $L(s)$ and $R(s)$, be defined according to equations \eqref{eq:singular-part} and \eqref{eq:regular-remainder}, respectively. 
Suppose that  $q_m(s)$ denotes  the degree-$m$ Chebyshev approximation of $R''(s)$
on $I_\delta$ as shown in equation \eqref{eq:Cheb-qm},  with $e_m$ denoting the corresponding error as defined in \eqref{eq:Cheb-error}.  Let us consider the  polynomial approximation $P_{m+2}(s)$ of $R(s)$ given by relation \eqref{eq:Pm-definition}, with $f_m$ representing the corresponding entropic approximation \eqref{eq:fm-definition}. 

Then
$
e_m\longrightarrow0
\quad\text{as }m\to\infty,
$
and the following properties hold.

\begin{enumerate}
\item
The regular part $R(s)$ is approximated in $C^2(I_\delta)$, that is, 
$
\|P_{m+2}-R\|_{C^2(I_\delta)}
\le C_\delta e_m,
$ holds, 
and 
\[
f_m\longrightarrow f
\qquad\text{in }C^2(I_\delta).
\]

\item Let $ \Phi_m(s)
=
f_m(s)-f_m(0)-\frac12s f_m'(s). $
For all sufficiently large $m$, there exists a coexistence sequence
$\{s_{NI}^{(m)}\}_{m\geq 1}$ satisfying
\[
\Phi_m(s_{NI}^{(m)})=0, \quad \text{and} \quad s_{NI}^{(m)}\longrightarrow s_{NI}, \,\, m\to\infty.
\]
\item
Let
\[
\kappa_m
=
\frac{3\theta_{NI}}
     {4s_{NI}^{(m)}}
f_m'(s_{NI}^{(m)}), \quad \text{and} \quad 
\psi_m(s,\theta)
=
\theta f_m(s)-\frac23\kappa_m s^2.
\]
Then
$\lim_{m\to\infty}
\kappa_m=\kappa.$ Furthermore, at $\theta=\theta_{NI}$, the coexistence conditions hold:
\[
\partial_s\psi_m(0,\theta_{NI})=0,
\quad
\partial_s\psi_m(s_{NI}^{(m)},\theta_{NI})=0, \quad  \psi_m(0,\theta_{NI})
=
\psi_m(s_{NI}^{(m)},\theta_{NI}).
\]

\item
The endpoint singularities are satisfied, that is, 
$
\psi_m(s,\theta)\longrightarrow+\infty
\quad\text{as}\qquad
s\to \{-\frac12, 1\}.
$
\end{enumerate}
\end{theorem}
\medskip
\noindent\textit{Proof.} The proof is given in Section~S1 of the Supplemental Material.\hfill$\square$

\begin{remark}
For the uniaxial Ball--Majumdar potential, the positive
coexistence value $s_{NI}$ is unique.  Moreover, if it is a 
simple simple root, that is, it satisfies
$
\Phi'(s_{NI})\neq0,
$
then the conclusion of the theorem is sharper.   For all sufficiently
large $m$, the zero $s_{NI}^{(m)}$ is unique in a fixed neighborhood
of $s_{NI}$ and
$
|s_{NI}^{(m)}-s_{NI}|
\le C e_m.
$
Consequently,
$
|\kappa_m-\kappa|
\le C e_m.
$
This provides stability as well as a quantitative error estimate, but it is not needed as an additional
hypothesis in the formulation of the approximation procedure.
\end{remark}

\section{Ericksen thermodynamic theory.} 

\subsection{Isothermal Ericksen model}
The pair of fields $\{\bn, s\}$ corresponding to the order tensor $Q$ in equation  \eqref{eq:Quniaxial} represent the independent variables of the isothermal Ericksen model \cite{Ericksen1991},  upon neglecting flow.    In this theory,  $s=0$ corresponds to the isotropic state (with $\bn$ undefined), $s=1$ the state of perfect alignment and  $s=-\frac{1}{2}$ an arrangement with the molecules located on a plane perpendicular to $\bn$. 
The isothermal free energy, parametrized by the temperature $\theta>0$,  takes the form 
\begin{align}
    I (s, \mathbf{n}, \theta) =& \int_\mathrm{\Omega}\rho W (s, \nabla s, \mathbf{n}, \nabla \mathbf{n}, \theta) \, d \mathbf{x}, \quad \text{with} \\ 
    W(s, \nabla s, \mathbf n, \nabla\bn, \theta)&= \frac{k}{2} |\nabla s|^2+ s^2 W_{\mathrm{OF}} (\bn, \nabla\bn)+ \psiB(s, \theta), \label{eq:Ericksen}
\end{align} 
with $W_{\mathrm{OF}}=\frac{k}{2}|\nabla\bn|^2$, in the so called one-constant case, 
that is, when the splay, twist and bend elastic constants are equal.  For most purposes, we will take $\psiB$ as in \eqref{energy_BM_uniaxial}, also allowing for  the  form \eqref{psi_unaxial}  in order to take  advantage of the polynomial Landau-de Gennes form for specific calculations.
Existence of minimizer of the energy \eqref{eq:Ericksen} has been studied by F. H. Lin \cite{Lin1991,LinPoon1994} and its gradient flow in \cite{CaldererGolovatyLinLiu2002}. Calculating the first variation with respect to the fields of the problem $\{s, \bn\}, $  we arrive at the Euler-Lagrange equations
\begin{equation}
    \Gelast:=\rho \big(\nabla\cdot\partderiv{W}{\nabla s}-\partderiv{W}{s}\big)=0, \quad \gbelast:= \rho\big(\nabla\cdot\partderiv{W}{\nabla\bn}-\partderiv{W}{\bn}\big)-\lambda\bn=0. \label{eq:E-equilibrium}
\end{equation}
The scalar and vector field $\Gelast$ and $\gbelast$ represent the generalized elastic forces, respectively, associated with the orientational fields $s$ and $\bn$, respectively.   Moreover, the first term in each expression corresponds to a force, in each case associated with a vector $\mathbf T$ and a tensor $\boldsymbol{\tau}$, respectively, 
\begin{equation} \label{eq:surface_forces}
    \mathbf T= \rho \partderiv{W}{\nabla s}, \quad \boldsymbol\tau= \rho\partderiv{W}{\nabla\bn}. \end{equation}

\subsection{Ericksen model with variable temperature}
We now  proceed to establishing a thermodynamic model for uniaxial nematic liquid crystals for the director, order parameter and temperature fields, 
\begin{equation} \bn=\bn(\bx, t), \, s=s(\bx, t), \, \theta=\theta(\bx, t), \quad \bn\cdot\bn=1,  \label{variables}\end{equation}
respectively. The fields $\{\dot \bn, \dot s\}$ are the generalized velocities, with $\nabla\theta$ also representing a dissipative field.
  The constitutive equations of the model are listed as
$$\Psi = \Psi(s, \theta, \nabla s, \bn, \nabla\bn), \, \,\eta=\eta(s, \theta, \nabla s, \bn, \nabla\bn),\,\, e=e(s, \theta, \nabla s, \bn, \nabla\bn), \,\, \bq=\bq(s, \theta, \nabla s, \bn, \nabla\bn, \nabla\theta),$$
denoting, the free energy density, the entropy and internal energy densities per unit mass, respectively,  with $\mathbf q$ representing the heat flux vector. They   are not all independent but satisfy the defining relation
\begin{equation}
    e= \Psi+\eta\theta, \label{eq:internal-energy}
\end{equation}
representing the energy contribution that can be transformed into mechanical work plus that lost in heat. 

The total energy of the system is defined as 
\begin{equation}
    \mathcal E=\int_\mathrm{\Omega}\rho e(\bx, t)\,d\bx. 
\end{equation}
The local form of the  law of balance of energy states that  the relation
\begin{equation} \label{eq:balance-energy}
    \frac{d}{dt}\int_{\mathcal U}\rho e(\bx, t)\,d\bx= \int_{\partial\mathcal U} [(\boldsymbol\tau\dot\bn +\mathbf T\dot s) - \bq]\cdot\boldsymbol\nu\, dS+ \int_{\mathcal U} \big(\gbext\cdot\dot \bn+ \Gext \dot s+r\big)\,d\bx,
\end{equation}
hold for every subset $\mathcal U\subseteq \mathrm{\Omega}$, 
where $\mathbf T$ and $\boldsymbol \tau$  are surface torques, $r$ denotes the volumetric heat supply, $\gbext$ and $\Gext$ external generalized forces. The expressions of $\mathbf T$ and $\boldsymbol \tau$ follow those in \eqref{eq:surface_forces} but  with $W$ replaced with $\Psi$, with the analogous expression for the generalized forces. 
The local form of the Second Law of Thermodynamics in the form of the Clausius-Duhem inequality establishes that
\begin{equation}
\frac{d}{dt}\int_{\mathcal U} \rho \eta \,d\bx\geq -\int_{\partial \mathcal U} \frac{1}{\theta} \bq\cdot\boldsymbol{\nu}\,dS+ \int_{\mathcal U}\frac{r}{\theta}\,d\bx, \label{eq:Clausius-Duhem}
\end{equation}b
holds for all $\mathcal U\subseteq \mathrm{\Omega}$.
Another source of energy dissipation, in addition to heat conduction, is due to viscous generalized forces (torques) $\{\hat \bg, \hat G\}$. 
The usual time-dependent theories of nematic postulate a linear dependence of the viscous generalized forces on the (generalized) velocities, 
that is,  $\hat \bg= -\gamma_1(s) \dot \bn$ and $\hat G=-\beta(s) \dot s, $ where $ \gamma_1(s)$ and $\beta(s)$ are known as rotational viscosities. 
We now postulate the laws of balance of generalized forces  (physically related to the balance of angular momentum) but neglecting the moment of inertia (with a contribution of the order of $10^{-18}$, with respect to the remaining terms).  These equations state that the sum of elastic,  (generalized) viscous and external forces, satisfy $\gbelast + \hat\bg+\gbext=0$ and $\Gelast + \hat G+\Gext=0$, with $\gbelast$ and $\Gelast$  as in \eqref{eq:E-equilibrium}. That is,  
\begin{align}
\gamma_1(s) \dot\bn=& \rho\big(\nabla\cdot\partderiv{\Psi}{\nabla\bn}-\partderiv{\Psi}{\bn}\big)-\lambda\bn+\gbext,  \label{eq:balance-n}\\
\beta(s)\dot s=&  \rho\big(\nabla\cdot\partderiv{\Psi}{\nabla s}-\partderiv{\Psi}{s}\big)+\Gext. \label{eq:balance-s}
\end{align}  
The local form of the equation of balance of energy \eqref{eq:balance-energy} gives

  \begin{equation}
    \rho \dot e= \nabla\cdot (\boldsymbol \tau\dot\bn + \mathbf T\dot s) - \nabla\cdot \bq+ \gbext\cdot\dot \bn+ \Gext \dot s+r  \label{eq:energy-local}
\end{equation}
 Substituting equations \eqref{eq:balance-n} and \eqref{eq:balance-s} into the previous one and simplifying the result, we get 
\begin{equation} {\label{eq:balance-energy-local-reduced}}
    \rho \dot e= (\gamma_1\dot\bn+ \rho\partderiv{\Psi}{\bn})\cdot\dot\bn +\boldsymbol{\tau}\cdot \nabla\dot\bn + (\beta\dot s+\rho \partderiv{\Psi}{s})\dot s+ \boldsymbol{T}\cdot\nabla\dot s-\nabla\cdot\bq. \end{equation}
We assume that the Fourier law of heat conduction
\begin{equation}
    \bq= -\kappa\nabla \theta, \quad \kappa= \kappa_\perp(s) \mathbf{I} +\kappa_a(s)(\bn\otimes\bn), \label{eq:thermal-conductivity}
\end{equation}
holds, where $\kappa $ represents the positive definite conductivity tensor for uniaxial nematic liquid crystals.  $\kappa_\perp>0$ and $\kappa_\|>0$ are scalar values of the conductivity along the directions perpendicular and parallel to $\bn$, respectively, and $0<\kappa_a= \kappa_\|-\kappa_\perp$. 
We now follow the customary process of deriving the heat equation from the balance of energy \eqref{eq:balance-energy-local-reduced} using the relation \eqref{eq:internal-energy} together with that following from the application of the Clausius-Duhem inequality, 
\begin{equation}
    \eta=-\partderiv{\Psi}{\theta}. \label{eq:entropy}
\end{equation}
This yields the  following:
\begin{align}
  -\theta\partderivtwo{\Psi}{\theta}\dot\theta  - \theta\partderivmix{\Psi}{\theta}{\bn}\cdot\dot\bn - \theta\partderivmix{\Psi}{\theta}{s}\dot s -\theta\partderivmix{\Psi}{\theta}{\nabla \bn}\cdot\nabla\dot \bn  -\theta\partderivmix{\Psi}{\theta}{\nabla s}\cdot\nabla\dot s
    = \frac{1}{\rho}\nabla\cdot(\kappa\nabla \theta). \label{eq:balance-e-long}
\end{align}
Following \cite{Ericksen1991}, we now take  $\Psi$ as
\begin{equation}
 \Psi(s, \nabla s, \bn, \nabla\bn, \theta)=   W(\nabla\bn, \bn, s, \theta)+H(s, \theta), \label{psi_modified}\end{equation}
with $W$ as in \eqref{eq:Ericksen}; $H(s, \theta)$, to be determined next,  has been added  so as to satisfy required  thermodynamics  properties, related to the positiveness of the specific heat capacity. 
 Substituting  \eqref{psi_modified}  into \eqref{eq:balance-e-long} gives
 \begin{align}
   -\rho \theta(\partderivtwo{\psiB}{\theta}+ \partderivtwo{h}{\theta}(\theta, s))\dot\theta 
    = \nabla\cdot(\kappa\nabla \theta) +\rho  \theta(\partderivmix{\psiB}{\theta}{s}+ \partderivmix{H}{\theta}{s})\dot s,   \label{eq:balance-e-short}
\end{align}
with $\psiB$ either as in \eqref{energy_BM_uniaxial} or \eqref{psi_unaxial}. We point out that the first two terms of the previous equation fit the format of the heat equation provided
\begin{equation}
  0<  c_\theta(s, \theta):=\partderiv{e}{\theta}=  -\theta(\partderivtwo{\psiB}{\theta} +\partderivtwo{H}{\theta}(\theta, s))= - \theta \partderivtwo{H}{\theta}(\theta, s). \label{c-theta} 
\end{equation}
Integrating  the last equation with respect to $\theta$ twice, yields
\begin{equation} \label{f-theta-modified}
    H(s, \theta) = \int_{\theta^*}^\theta\int_{\theta^*}^r -\frac{{c_\theta}(s, u)}{u}\,du\,dr.
\end{equation} 
In the case that  $c_\theta(s, \cdot)$ is independent of $s$; its graph is shown in
Figure \ref{fig:double_well_transition}(a). On the other hand, if it  is independent of $\theta$,  then 
\begin{equation} \label{f-theta}
    \theta \partderivtwo{H}{\theta}=-c_\theta(s) \implies h(\theta, s) = -{c_\theta}(s) \theta \ln\theta 
\end{equation} 
Note that, in either case, the positivity of $c_\theta$ implies concavity of $h(s, \cdot)$.

In summary, the system of governing equations is given by 
\begin{align} \label{governing-system}
\begin{cases}
     \gamma_1(s) \dot\bn = \rho\big(\nabla\cdot\partderiv{W}{\nabla\bn}-\partderiv{W}{\bn}\big)-\lambda\bn \\
    \beta(s)\dot s = \rho\big(\nabla\cdot\partderiv{W}{\nabla s}-\partderiv{W}{s}-\partderiv{H}{s}\big) \\
    \rho c_\theta \dot\theta = \nabla\cdot(\kappa\nabla \theta) + \rho \theta(\partderivmix{\psiB}{\theta}{s}+\partderivmix{H}{\theta}{s})\dot s. 
\end{cases}
\end{align}
We now set the domain of the problem to be the experimental  cylinder
\begin{equation}
    \mathrm{\Omega}=\{\bx=(r, \phi, z); 0<r<\frac{l}{2}, \, \,0<\phi<2\pi, \,\,0<z<d\}, \quad \partial\Omega=\Gamma\cup\Gamma_-\cup\Gamma_+, 
\end{equation}
represented in figure \ref{fig:cylinder_domain}.  We let $\Gamma$ denote the lateral surface and $\Gamma_+$ and $\Gamma_-$, the basis at the top and bottom of the cylinder, respectively. 

\subsubsection{Boundary and Initial Conditions}
The  boundary conditions for the optic variables $\bn, s$  are prescribed as follows:  
\begin{align}
   &\frac{\partial \bn}{\partial{\boldsymbol \nu}} =0, \, \text{on}\, \Gamma_-\cup\Gamma_+, \quad \bn={\mathbf e}_z, \, \text{on}\, \Gamma, \label{eq:bc-n}\\
    & \frac{\partial s}{\partial{\boldsymbol \nu}} =0\, \, \text{on}\, \partial{\Omega}, \label{eq:bc-s} 
\end{align}
where $\boldsymbol \nu$ is the unit outer normal to the boundary. This indicates that the liquid 
crystal is anchored along the direction parallel to the axis on the lateral boundary, $\Gamma, $ of ${\Omega}$ whereas it is left free on the bases $\Gamma_-$ and $\Gamma_+$ (the latter correspond to the natural boundary conditions stemming from minimization of the energy). 

Let us now set up the thermal boundary conditions. 
In laboratory experiments, the liquid crystal sample is placed in contact with a temperature-controlled stage (e.g., a Peltier device), whose temperature is externally prescribed and varied in time. The apparatus controls the temperature of the surrounding thermal reservoir rather than the temperature of the sample itself. As a result, heat exchange between the sample and the environment occurs through a finite thermal contact resistance. This process is modeled by a Newton-Robin boundary conditions,
\begin{gather}
 \bq\cdot \boldsymbol{\nu} :=\zeta (\theta-\theta_{\text{b}}(t)), \quad \text{on}\,\, \Gamma_-\label{eq:bc-Robin}\\
\theta_\text{b}(t)=\theta^+ -h t,  \quad h>0, \,\, \text{for}\, \, t>0, \label{eq:bc-cooling}
\end{gather}
where $\theta^+>\theta_{IN}$ is a reference temperature and 
$\zeta>0$ the effective convective heat transfer coefficient of the material in the cooling bath. 
\begin{equation}
    \bq\cdot \boldsymbol{\nu} = 0, \quad \text{on}\,\, \Gamma_+\cup\Gamma. \label{eq:bc-insulated}
\end{equation}
\begin{remark} The  
Dirichlet boundary condition 
\, $
    \theta(-\frac{d}{2}, t)= \theta^+ -h t,  $ 
alternate to \eqref{eq:bc-Robin}, corresponds to the perfect transmission limit $\zeta\to\infty$. 
\end{remark}
We prescribe  initial data 
\begin{equation}
  s^{{\tiny \rm in}}(\bx)=  s(\bx, 0), \, \, \theta^{{\tiny \rm in}}(\bx)= \theta(\bx, 0),  \,\,  \bn(\bx, 0)= \mathbf e_z. \label{eq:initial}
\end{equation}
Experimental values of the parameters are given in Table \ref{table:1}.
\begin{remark}
    The function 
 $  \,  V(t)=\int_{\Omega}\rho\big(e(\bx,t)-\theta_\infty \eta(\bx, t)\big)\,d\bx, \, $
where $\theta_\infty>0$ is a reference temperature,  is a Lyapunov function of the system provided $h=0$ in \eqref{eq:bc-cooling}, that is, the bath remains at constant temperature and no cooling takes place.   
\end{remark}

\subsection{Initial-Boundary Value Problem.} \label{sec:iBVP}

\begin{proposition} Suppose that the initial and boundary conditions \eqref{eq:bc-n} and \eqref{eq:initial} hold. Then the system \eqref{governing-system} admits the solution $\bn(r, \theta, z, t) =\mathbf e_z$ in $\Omega$ for $t\geq 0.$
\end{proposition}
The argument behind the proposition is that one can maintain a vertical  director field  on the boundary by  rubbing the surface along such a direction. This does not contradict the fact that, at high temperature the liquid crystal is still isotropic, with $s=0$ indicating that the actual molecules align poorly with  $\mathbf e_z$. 

Consequently, the governing equation \eqref{eq:balance-n} admits the solution $\bn=\mathbf e_z$, upon the appropriate choice of the Lagrange multiplier $\lambda$. 
Note that the free energy reduces to
 \begin{align}
  \Psi(s, \nabla s, \bn, \nabla\bn, \theta)= \frac{K}{2} |\nabla s|^2+ \psiB(s, \theta)+ H(s,\theta), \quad    
  \label{eq:LdGs}
\end{align}
with $\psiB$  and $h(s,\theta) $ as in \eqref{psi_unaxial} and \eqref{f-theta-modified}, respectively. 
Furthermore, with  $c_\theta, \kappa, \beta >0$, the previous equations form a   nonlinear parabolic system, 
 \begin{align}{\rho } c_\theta\dot\theta
    = &\nabla\cdot(\kappa\nabla \theta) -(\partderiv{\psiB}{s}+\partderiv{H}{s}- \theta\partderivmix{\psiB}{\theta}{s}-\theta \partderivmix{H}{s}{\theta})\dot s 
\label{eq:balance-e-simplified}\\
    \beta(s)\dot s=& \big(K\triangle s-\partderiv{\psiB}{s}(s,\theta)\big). \label{eq:balance-s-simplified}
\end{align}

\begin{remark}
In the case that $c_\theta(\theta)$ is independent of $s$, equation  \eqref{eq:balance-e-simplified} should be replaced with 
\begin{equation}\label{eq:balance-e-indep_of_s} {\rho } c_\theta\dot\theta
    = \nabla\cdot(\kappa\nabla \theta) -(\partderiv{\psiB}{s}- \theta\partderivmix{\psiB}{\theta}{s})\dot s.  
\end{equation}
\end{remark}
 and varied in time. The apparatus controls the temperature of the surrounding thermal reservoir
rather than
\begin{table}[h!]
\centering
\begin{tabular}{||c c c c ||} 
 \hline
 Quantity & Name & Units & Dimensional values \\ [0.5ex] 
 \hline\hline
 $d$ & Domain height & m & $3.44\times10^{-4}, 10^{-2}$ \\
 $l$ &  Domain diameter & m & $5.4\times10^{-3}$ \\
 vol & Domain Volume & $\text{m}^3$ &7.874 $ \times 10^{-9}$ \\
 $k$ & Frank constant & N & ${10}^{-12}$ \\
 $c$ & Bulk energy quartic coefficient & J/$\text{m}^3$ & $0.35\times10^6$ \\
 $b$ & Bulk energy cubic coefficient & J/$\text{m}^3$ & $-0.64\times10^6$\\ 
 $\alpha$ & Bulk energy quadratic coefficient & J/(${\text m}^3\text{K})$ & $b^2/(351c)=3.3341\times10^3$ \\
 $\rho$ & Mass density & kg/${\text m}^3$ & $10^3$\\
 $c_\theta$ & Specific heat capacity & J/(kgK)& { 2300} \\
 $ \beta$ & Orientational viscosity & Pa $\cdot$ s & {0.1}\\ 
 $\kappa$ & Thermal conductivity & W/(${\text{m}. \text{K}}$) & 0.8  \\
 $\frac{\kappa}{\rho c_\theta}$ & Thermal diffusivity & $\frac{\text{m}^2}{s}$ & {3.478 $\times 10^{-7}$} \\
 $\theta_{\textrm{m}}$ & Melting temperature & K & 439 \cite{Teerakapibal2018} \\
 $\theta_{\textrm{NI}}$ & Isotropic-nematic transition temperature & K & 363 \cite{Teerakapibal2018} \\ 
 $\theta_{\textrm{N/Sm}}$ & Nematic-smectic transition temperature & K & 347 \cite{Teerakapibal2018} \\ 
 $\theta_{\textrm{g}}$ & Glass transition temperature & K & 328 \cite{Teerakapibal2018} \\
$h$ & Cooling rate & K/s & $0.0167$\\
  $t_\circ$ & Reference time ($K/h$) & s & 59.86 \\
  $\zeta$ & Convective heat transfer coefficient & W/($\textrm{m}^2$.K) & 100 \cite{carino2013evaluation} \\
 \hline
\end{tabular}
\caption{Values of parameters for the ITZ liquid crystal. Note that  $t_0$ corresponds to the time that it takes the bath water to lower its temperature by 1 second. There is also another time scale associated with the dynamics of the liquid crystal, this calculated in terms of the generalized viscosity. This attests to the multiscale (time and also space) nature of the problem. We have used the value of  $c_\theta$ near the isotropic nematic transition.
The value of $\beta$ is taken from the liquid crystal literature \cite{Baron2001}, \cite{GruzdenkoDierking2023}, \cite{Khazimullin2011}, \cite{SchadZeller1982}.}
\label{table:1}
\end{table}
\begin{remark}

We identify the kinetic coefficient $\beta$ in the evolution equation for the scalar
order parameter with the rotational viscosity
$\gamma_1$ either in the Ericksen-Leslie or the Beris-Edward model for (uniaxial) $Q$.
 For low-molecular-weight
thermotropic nematic liquid crystals, measured values typically lie in
the range
$\,
10^{-2}\ {\rm Pa\,s}
\;\lesssim\;
\gamma_1
\;\lesssim\;
1\ {\rm Pa\,s},$
with values increasing as the nematic--isotropic transition is
approached. Consequently, simulations based on
$
\beta=\gamma_1=0.05,\;0.1,\;0.5,\;1\ {\rm Pa\,s}$
span the range commonly reported for thermotropic nematics. Although such a coefficient for ITZ is presently unknown, since it
is a highly viscous  material, its effective
rotational viscosity may differ substantially from that of classical
nematics.
\end{remark}
\begin{remark}
The coefficient $c_\theta$ in \eqref{c-theta} represents the specific heat capacity with the orientational state
held fixed. Experimental calorimetric measurements are generally reported
as the constant-pressure specific heat $c_p$. Since ITZ is treated
here as a nearly incompressible condensed material, with negligible
mechanical deformation, the difference between $c_p$ and $c_{\theta}$ is small.
We therefore identify
$c_\theta\simeq c_p.
 $
High-resolution adiabatic scanning calorimetry measurements give a heat-capacity of approximately
$
c_p=2.2\text{--}2.4\ {\rm J\,g^{-1}K^{-1}}
$
in the nematic and isotropic temperature range. We consequently use the
representative value
$
c_\theta=2.3\times10^3\ {\rm J\,kg^{-1}K^{-1}}.
$
The pronounced pretransitional increase in the measured equilibrium heat
capacity near the isotropic--nematic transition is not included in this
baseline value. 
\end{remark}

\section {Existence and uniqueness of solution}
 Also, add the double well profile.  Change notation $W_0$, not used anymore

In this section, we study the well-posedness of the initial-boundary
value problem corresponding to equations
\eqref{eq:LdGs}, \eqref{eq:balance-e-simplified},
\eqref{eq:balance-s-simplified}, and \eqref{eq:initial}, with the
boundary conditions described below.

The physical sample occupies a finite cylindrical domain. In order to
avoid the regularity issues associated with the circular edges where
the lateral boundary meets the top and bottom bases, we work with a
smooth regularization $\Omega$ of the cylinder, obtained by rounding
these edges in arbitrarily small neighborhoods. The boundary
conditions are properly regularized as well in these neighborhoods.

More precisely, for the thermal boundary condition we choose
$\chi\in C^\infty(\partial\Omega)$, with $0\leq\chi\leq1$, such that
$\chi=1$ on the portion of the boundary corresponding to the bottom
base $\Gamma_-$, away from the rounded edge, and $\chi=0$ on the insulated
portion, $\Gamma\cup \Gamma_+$ also away from the rounded edge. On the regularized boundary we
prescribe

\begin{align}
\partial_\nu s &=0
&&\text{on }\partial\Omega,
\label{eq:bc-s-reg}\\
-\kappa_T\partial_\nu\theta
&=\zeta\chi(x)\bigl(\theta-\theta_b(t)\bigr)
&&\text{on }\partial\Omega.
\label{eq:bc-theta-reg}
\end{align}

Thus the physical Robin and Neumann conditions are recovered away
from the rounded neighborhoods, while the transition between them is
smooth. The regularized problem therefore has a smooth boundary and
a single smooth Robin-type boundary relation.
\bigskip
\begin{theorem}[Local existence, uniqueness, separation, and continuation]
\label{thm:BM}

Let $\Omega\subset\mathbb R^3$ be the smooth regularization of the
cylindrical domain described above, and let
$
p>3.$
Assume
$
\rho c_\theta>0,\, \,
\kappa_T>0,\,\,
\beta>0,\,\,
K>0,\,\,
\zeta>0.
$
Assume that the prescribed bath temperature is spatially homogeneous,
with
\[
\theta_b\in W^{1,p}_{\rm loc}([0,\infty)),
\qquad
\theta_b(t)\ge\vartheta_b>0.
\]
Let us assume that the initial data satisfy
$
\theta^{{\tiny \rm in}},s^{{\tiny \rm in}}\in W^{2-2/p,p}(\Omega),
$
and  that, for some $\vartheta_0>0$ and
$\varepsilon_0>0$,
\[
\theta^{\rm in}(x)\ge\vartheta_0,
\qquad
-\frac12+\varepsilon_0
\le s^{{\tiny \rm in}}(x)\le
1-\varepsilon_0,
\qquad x\in\overline\Omega.
\]
Moreover the initial data are assumed to satisfy the compatibility conditions
\[
\partial_\nu s^{\rm in}=0,
\qquad
-\kappa_T\partial_\nu\theta^{\rm in}
=
\zeta\chi(x)
\bigl(\theta^{\rm in}-\theta_b(0)\bigr)
\qquad\text{on }\partial\Omega.
\]
Then there exists $T_*>0$ and a unique strong solution
$(\theta,s)$ of equations \eqref{eq:balance-e-simplified} and
\eqref{eq:balance-s-simplified}, with the  bulk
potential \eqref{eq:LdGs}, subject to the regularized boundary conditions
\eqref{eq:bc-s-reg}--\eqref{eq:bc-theta-reg} and the initial condition
\eqref{eq:initial}, 
 such that
\[
\theta,s\in
W^{1,p}(0,T_*;L^p(\Omega))
\cap
L^p(0,T_*;W^{2,p}(\Omega)).
\]
Moreover,
\[
\theta,s\in
C\bigl([0,T_*];W^{2-2/p,p}(\Omega)\bigr)\cap C(\overline\Omega\times[0,T_*]).
\]
Upon decreasing $T_*$, if necessary, there exist constants
$
0<\vartheta_*<\vartheta_0,
\, \, 
0<\varepsilon_*<\varepsilon_0
$
such that
\[
\theta(x,t)\ge\vartheta_*,
\qquad
-\frac12+\varepsilon_*
\le s(x,t)\le
1-\varepsilon_*
\]
for all $(x,t)\in\overline\Omega\times[0,T_*]$.
Furthermore the solution depends continuously on the initial and prescribed
boundary data in the corresponding maximal-regularity topology.

The solution admits a unique maximal continuation to an interval
$[0,T_{\max})$, where
\[
T_*\le T_{\max}\le\infty.
\]
If $T_{\max}<\infty$,  then at least one of the following
alternatives holds:
\begin{align*}
&\limsup_{t\uparrow T_{\max}}
\left(
\|\theta(t)\|_{W^{2-2/p,p}(\Omega)}
+
\|s(t)\|_{W^{2-2/p,p}(\Omega)}
\right)
=\infty, \quad \text {or} \\
&\liminf_{t\uparrow T_{\max}}
\min_{x\in\overline\Omega}\theta(x,t)=0, \quad \text{or}\quad 
\liminf_{t\uparrow T_{\max}}
\operatorname{dist}
\left(
s(\overline\Omega,t),
\left\{-\frac12,1\right\}
\right)=0.
\end{align*}
If
$
\inf_{(x,t)\in\Omega\times[0,T_{\max})}\theta(x,t)>0,
$
then the last alternative is excluded.
\end{theorem}

\begin{proof}
Let
\[
g(s):=f'(s),\qquad -\frac12<s<1.
\]
Then $g\in C^\infty((-1/2,1))$, with the properties
\[
g(s)\to-\infty\quad\text{as }s\downarrow-\frac12,
\qquad
g(s)\to+\infty\quad\text{as }s\uparrow1.
\]
We give the main steps of the proof; complete details are provided in
Section~S2 of the Supplementary Material.

\medskip
\noindent\textbf{Step 1. Linear setting and truncation.}
For $T>0$ set
\[
E_p(T)=W^{1,p}(0,T;L^p(\Omega))
       \cap L^p(0,T;W^{2,p}(\Omega)).
\]
Maximal $L^p$-parabolic regularity for the Neumann and Robin
initial and boundary-value problems  of the heat operator gives the corresponding linear
estimates (\cite{amann1995linear}, \cite{Pruss2002}), while, since $p>3$,
\[
E_p(T)\hookrightarrow
C([0,T];W^{2-2/p,p}(\Omega))
\hookrightarrow C(\overline\Omega\times[0,T]).
\]
Choose $0<\varepsilon<\varepsilon_0$ and replace $g$ temporarily by
$g_\varepsilon\in C^\infty(\mathbb R)$ agreeing with $g$ on
$[-1/2+\varepsilon,1-\varepsilon]$.

\medskip
\noindent\textbf{Step 2. Fixed-point construction.}
Fix $T_0>0$ and let
$(\theta^0,s^0)\in E_p(T_0)^2$ solve the homogeneous linear
parabolic equations with the prescribed initial and boundary data.
For $0<T\le T_0$, set
\[
{}_0E_p(T)=\{v\in E_p(T):v(\cdot,0)=0\}
\]
and
\[
\mathcal B_R(T)=
\left\{
(\bar\theta,\bar s)=(\theta^0,s^0)+(u,v):
u,v\in{}_0E_p(T),\
\|u\|_{E_p(T)}+\|v\|_{E_p(T)}\le R
\right\}.
\]
For $(\bar\theta,\bar s)\in\mathcal B_R(T)$, solve successively
\begin{align}
\beta s_t-K\Delta s
&=-\bar\theta g_\varepsilon(\bar s)
  +\frac43\kappa_B\bar s,                                      \label{eq:fp-s}\\
\rho c_\theta\theta_t-\kappa_T\Delta\theta
&=\frac43\kappa_B\bar s\,s_t,                                  \label{eq:fp-theta}
\end{align}
with the prescribed initial and boundary data, and define
\[
\mathcal T(\bar\theta,\bar s)=(\theta,s).
\]
Maximal regularity makes $\mathcal T$ well defined.

The main short-time estimate is that, for some $\eta>0$,
\begin{equation}
\label{eq:small-time}
\|v\|_{L^\infty(0,T;C(\overline\Omega))}
\le CT^\eta\|v\|_{E_p(T)},
\qquad v\in{}_0E_p(T),
\end{equation}
with $C$ independent of $0<T\le T_0$.  Indeed, choosing $\alpha $ and $\eta$ such that
\[
\frac{3}{2p}<\alpha<1-\frac1p,
\qquad
0<\eta<1-\alpha-\frac1p,
\]
the mixed-derivative theorem gives
\[
E_p(T)\hookrightarrow
W^{1-\alpha,p}(0,T;W^{2\alpha,p}(\Omega))
\hookrightarrow C^{0,\eta}([0,T];C(\overline\Omega)),
\]
and \eqref{eq:small-time} follows from the zero initial trace.

Now, writing $w=s-s^0$ and $z=\theta-\theta^0$, maximal regularity gives
\[
\|w\|_{E_p(T)}\le C_R T^{1/p},
\qquad
\|z\|_{E_p(T)}
\le C_R\left(
\|(s^0)_t\|_{L^p(\Omega\times(0,T))}+T^{1/p}
\right).
\]
Hence $\mathcal T$ maps $\mathcal B_R(T)$ into itself for $T$
sufficiently small.

If $U_i=(\bar\theta_i,\bar s_i)\in\mathcal B_R(T)$, subtraction of
the two problems \eqref{eq:fp-s}--\eqref{eq:fp-theta}, together with
\eqref{eq:small-time}, the local Lipschitz continuity of
$g_\varepsilon$, and maximal regularity, yields
\[
\|\mathcal T(U_1)-\mathcal T(U_2)\|_{E_p(T)^2}
\le
CT^\eta\|U_1-U_2\|_{E_p(T)^2}.
\]
Thus, after decreasing $T$ so that $CT^\eta<1$, Banach's fixed-point
theorem gives a unique strong solution of the truncated problem.

\medskip
\noindent\textbf{Step 3. Removal of the truncation and separation.}
The embedding $E_p(T)\hookrightarrow
C(\overline\Omega\times[0,T])$ and the strict positivity and
separation of the initial data imply, for sufficiently small $T$,
\[
\theta\ge\vartheta_*>0,\qquad
-\frac12+\varepsilon_*\le s\le1-\varepsilon_*,
\]
with $\varepsilon_*>\varepsilon$.  Hence $g_\varepsilon(s)=g(s)$,
so the truncation is inactive and the solution solves the original
system.

More generally, let us suppose that a strong solution exists on $[0,T]$ and
$\theta\ge\underline\vartheta>0$.  Writing
\[
\beta s_t-K\Delta s+F(\theta,s)=0,
\qquad
F(\theta,s)=\theta g(s)-\frac43\kappa_Bs,
\]
the singular limits of $g$ imply that, for some $\delta>0$,
\[
F(\theta,s)\le0
\quad\text{near }s=-\frac12,
\qquad
F(\theta,s)\ge0
\quad\text{near }s=1,
\]
uniformly for $\theta\ge\underline\vartheta$.  Multiplying the latter equation by the  test functions 
$-(a-s)_+$ and $(s-b)_+$, where
$a=-1/2+\delta$ and $b=1-\delta$, and integrating  gives
\[
\frac{\beta}{2}\frac{d}{dt}\|(a-s)_+\|_2^2
 +K\|\nabla(a-s)_+\|_2^2\le0,
\]
and
\[
\frac{\beta}{2}\frac{d}{dt}\|(s-b)_+\|_2^2
 +K\|\nabla(s-b)_+\|_2^2\le0.
\]
Since the positive parts vanish initially,
\[
-\frac12+\delta\le s(x,t)\le1-\delta.
\]

\medskip
\noindent\textbf{Step 4. Uniqueness and continuation.} We note that, 
on compact subsets of
$(0,\infty)\times(-1/2,1)$,  the nonlinearities are locally Lipschitz.
Applying the same maximal-regularity and zero-trace estimates to the
difference of two solutions gives uniqueness on a short interval,
and iteration yields uniqueness on the common interval of existence.
Separating the linear solution carrying differences in the initial
data gives continuous dependence.

Finally, if $[0,T_{\max})$ is the maximal interval of existence and
the trace norms remain bounded while $(\theta,s)$ stays in a compact
subset of $(0,\infty)\times(-1/2,1)$, the preceding local argument
can be restarted with uniform constants.  Since $p>3$, the boundary
compatibility conditions pass to each time slice.  The solution can
therefore be extended beyond $T_{\max}$, a contradiction.  This
proves the stated continuation alternatives; the separation argument
above excludes approach to $s=-1/2$ or $s=1$ whenever the temperature
remains uniformly positive.
\end{proof}

\begin{remark}[The cylindrical domain]
\label{rem:cylinder}
The regularization of the cylindrical edges is introduced only to place
the well-posedness argument within the standard parabolic theory for
smooth domains.  The original physical domain is a finite cylinder and
is therefore piecewise smooth, with circular edges at the intersections
of the lateral surface with the top and bottom bases.  The problem may
also be formulated directly on this domain by using the parabolic
regularity theory appropriate to domains with edges, together with the
corresponding compatibility conditions on the initial and boundary data
and, when required, function spaces adapted to the edge singularities.

This modification concerns only the linear regularity theory used in
the local existence argument.  Once the corresponding parabolic
estimates are available, the nonlinear part of the proof is unchanged:
in particular, the fixed-point argument and the separation property
generated by the Ball--Majumdar singularity apply in the same way.
We do not pursue the technical regularity theory for the nonsmooth
cylindrical domain here.
\end{remark}
\begin{remark}[Higher regularity]
Once the solution is strictly separated from the singular values
$s=-1/2$ and $s=1$, the bulk potential and all of its
derivatives are bounded on the range attained by $s$. Consequently,
standard parabolic regularity and bootstrapping may be applied.
Thus, if the initial and boundary data have additional regularity
and satisfy the corresponding compatibility conditions, the solution
inherits the associated higher parabolic regularity. We do not pursue
these higher-regularity statements here.
\end{remark}

\begin{remark}[Separation and continuation]
The asymptotic properties of the singular potential strengthen the continuation criterion
in Theorem~\ref{thm:BM}. As shown in the proof, on any interval on
which
\[
\theta(x,t)\ge\underline{\vartheta}>0,
\]
the scalar order parameter remains uniformly separated from the
singular values $-1/2$ and $1$. Hence the third alternative in the
continuation criterion cannot occur while the temperature remains
uniformly positive. In particular, if the Sobolev norm in the first
alternative remains bounded, finite-time breakdown can occur only
through loss of a uniform positive lower bound for the temperature.
\end{remark}

\begin{corollary}[Persistence of one-dimensional symmetry]
\label{cor:one-dimensional}
Suppose that the domain is the cylindrical sample described in
Section~D and that the boundary data are uniform on each horizontal
cross-section. Assume, in addition, that the initial data depend only
on the axial variable,
\[
  \theta^{\rm in}(x)=\Theta^{\rm in}(z),
  \qquad
  s^{\rm in}(x)=S^{\rm in}(z).
\]
Then, on its interval of existence, the unique strong solution of the
initial-boundary value problem depends only on $z$ and $t$:
\[
  \theta(x,t)=\Theta(z,t),
  \qquad
  s(x,t)=S(z,t).
\]
Consequently, $(\Theta,S)$ is the unique strong solution of the
corresponding one-dimensional initial-boundary value problem.
\end{corollary}

\begin{proof}
The one-dimensional problem obtained by restricting the governing
equations to fields depending only on $z$ admits, upon extension
constant in the transverse variables, a strong solution of the
three-dimensional problem with the prescribed initial and boundary
data. Uniqueness in Theorem~\ref{thm:BM} therefore implies that the
three-dimensional solution coincides with this extension throughout
their common interval of existence.
\end{proof}
The preceding corollary justifies the one-dimensional reduction used in the remainder of the paper. Since the experimental initial and boundary data are uniform across the horizontal cross-sections of the sample, we henceforth consider solutions depending only on the axial coordinate \(z\) and time \(t\).

\section{Scaling, non-dimensionalization and asymptotical analysis}
\label{sec:scaling}

We now study the motion of the front separating the nematic from the isotropic region as the cooling of the bath is underway. For this, we set the equations into dimensionless form and carry out their asymptotic analysis. This will allow us to explore the different regimes according to the thickness of the cylinder. It turns out that the thickness parameter also determines the time scale of the front motion. We determine a critical thickness below which the front evolves in the cooling time scale of the bath, whereas above it, the dynamics follows that of the heat diffusion in the material.

Let $t_0$ and $\theta^*$ denote the typical time and temperature scales, respectively.  We define the dimensionless variables
\begin{equation}
    \bar{x} = \frac{x}{l}, \quad \bar{y} = \frac{y}{l}, \quad \bar{z} = \frac{z}{d}, \quad \bar{t} = \frac{t}{t_0}, \quad \bar{\theta} = \frac{\theta}{\theta^*}, \quad \bar{h} = \frac{h t_0}{\theta^*},
    \label{eq:nondim_vars}
\end{equation}
where $l$ and $d$ represent the diameter and height of the cylindrical domain, respectively. 
Using these variables, the governing system of equations transforms into the following dimensionless form:
\begin{align}
    B_1 \frac{\partial\bar{\theta}}{\partial\bar{t}} &= \frac{\partial^2\bar{\theta}}{\partial\bar{z}^2} - \left[B_2 s^3 + B_3 s^2 - B_4 s\right] \frac{\partial s}{\partial\bar{t}}, \label{eq:temp_nondim} \\
    C_1 \frac{\partial s}{\partial\bar{t}} &= C_2 \frac{\partial^2 s}{\partial\bar{z}^2} - \left[(\bar{\theta}-1)s + C_3 s^2 + C_4 s^3\right]. \label{eq:order_nondim}
\end{align}
The boundary conditions representing Robin convective cooling at the bottom ($\bar{z}=0$) and thermal insulation at the top ($\bar{z}=1$), together with the natural conditions on $s$ become,
\begin{align}
  &  \left( \pdv{\bar{\theta}}{\bar{z}} \right)_{\bar{z}=0}=\textrm{Bi}(\bar{\theta}-\bar{\theta}^+ + \bar{h}\bar{t}), \label{eq:bc-theta_neumann} \\
&    \left( \pdv{\bar{\theta}}{\bar{z}} \right)_{\bar{z}=1}=0, \quad
 \left( \pdv{s}{\bar{z}} \right)_{\bar{z}=0}=0=\left( \pdv{s}{\bar{z}} \right)_{\bar{z}=1}. \label{eq:bc-s_neumann}
\end{align}
Similarly, the initial conditions are given by the following: 
\begin{equation}
    s(\bar{z},0)=\epsilon, 0<\epsilon <<1, \quad \bar{\theta}(\bar{z}, 0)=\bar{\theta}^+, \epsilon=\sqrt{B_1}. \label{eq:initial-s_theta} 
\end{equation}
 
The dimensionless parameter groups  are of the form
\begin{align} \label{parameter_groups}
&B_1 = \frac{\rho c_\theta d^2}{\kappa_{\parallel} \tau}, \,\, B_2 = \frac{8cd^2}{9 \theta^* \kappa_\parallel \tau}, \,\, B_3 = \frac{4bd^2}{9 \theta^* \kappa_\parallel \tau}, \,\, B_4 = \frac{4 \alpha d^2}{3 \kappa_\parallel \tau}, \,\,  \textrm{Bi}=\frac{\zeta d}{\kappa_\parallel}, \\  &C_1=\frac{3 \beta}{4\tau \alpha \theta^*}, \,\, C_2 = \frac{3K}{4 \alpha d^2 \theta^*}, \,\, C_3 = \frac{b}{3 \alpha \theta^*}, \,\, C_4 = \frac{2c}{3 \alpha \theta^*}, \,\, \bar h= \frac{h\tau}{\theta^*}, \end{align}
where $\zeta$ and $\text{Bi}$ denote the convective heat transfer coefficient and the Biot number, respectively. The latter represents the ratio of the convective heat transfer at the bottom wall of the cylindrical domain in contact with the cooling bath to the conductive heat transfer within the cylindrical domain of the liquid crystal. 

Next, we examine the time scales of the model and derive the critical thickness that separates the corresponding regimes. 
We  denote 
\begin{equation}
    \tau:=t_0= \frac{1^\circ \text{K}}{h}, \qquad \tau:=\tau_D = \frac{\rho c_\theta d^2}{\kappa_{||}}, \qquad B_1=\frac{\tau_D}{t_\circ}= \frac{\rho c_\theta d^2}{\kappa_{||} t_0}\label{time_scales}
\end{equation}
the time scales corresponding to the bath cooling and the diffusion rates, and their ratio, respectively, the latter governing the phase transition dynamics.   Since  $B_1 \propto d^2$, the choice of the appropriate time scale depends strictly on the physical dimensions of the cylinder. The transition between the cooling and the  diffusion governed  regimes occurs when $B_1=1$, that is
\begin{equation}
    d_{\text{crit}} = \sqrt{\frac{\kappa_{||} t_0}{\rho c_\theta}} = \sqrt{\alpha_T t_0},
    \label{eq:d_crit}
\end{equation}
where $\alpha_T \approx 3.478 \times 10^{-7}\text{ m}^2/\text{s}$ is the thermal diffusivity of the ITZ material, and $t_0 = 59.86\text{ s}$, for the baseline cooling rate given in Table 1, in which case, the critical height threshold is $d_{\text{crit}} \approx 4.56\text{ mm}$. 

This threshold distinguishes  two distinct regimes:
\begin{enumerate}
    \item \textbf{Bath-cooling ($d \ll d_{\text{crit}}$):} For thin cylinders, such as, $d = 0.344\text{ mm}$), $B_1 \ll 1$. Heat diffuses in the bulk almost instantaneously relative to the bath cooling, validating the $t_0$ time scale.
    \item \textbf{Bulk-Diffusion ($d > d_{\text{crit}}$):} For thick laboratory or industrial samples (e.g., $d \ge 1.0\text{ cm}$), $B_1 > 1$. The time scale of the front motion is dictated by internal heat conduction.
\end{enumerate}

\subsection{Bath-cooling time scale and asymptotic analysis}
For the parameter values in the table, with $d=0.344 \text{mm}$ and $t_\circ= 59.86$, the parameter set in \eqref{parameter_groups} yields
\begin{align}
B_1 &= 5.683\times10^{-3}, \quad
B_2 = 2.196\times10^{-6}, \quad 
B_3 = -2.008\times10^{-6}, \quad 
B_4 = 1.098\times10^{-5}, \quad \text{Bi}= 0.043\nonumber \\
C_1 &= 1.073\times10^{-9}, \quad
C_2 = 5.431\times10^{-12}, \quad 
C_3 = -0.1828, \quad  
C_4 = 0.1999, \quad \bar h= 2.855\times 10^{-3}.
 \label{parameter_set1}
\end{align}

 As a result of the parameter gap, we see that the following time and length scales will appear in the solutions of the system: 
\begin{equation}
    \{\bar t, \frac{\bar t}{B_1}, \frac{\bar t}{C_1}, \bar z, \frac{\bar z}{\sqrt{C_2}}\}.
\end{equation}
Let us first consider equation \eqref{eq:order_nondim} and note that the solution can be formally written as 
$$s(\bar z, \bar t, \frac{\bar t}{C_1}) = S(\bar z, \bar t; C_1, C_2) + S_I(\bar z, \frac{\bar t}{C_1}), $$ with $S_I$ denoting the initial layer correction to the outer solution. It is easy to check that there is no boundary layer contribution in the scale $\frac{\bar z}{\sqrt{C_2}}$ due to the homogeneous Neumann boundary conditions.
The outer solution $S$ belongs to one of two branches of the cubic equation
\begin{equation}
(\bar{\theta}-1)S+C_3 S^2 + C_4 S^3 = 0: \quad S=0 \,\,\text{and}\, \, S =  \frac{-C_3+\sqrt{C_3^2-4(\bar{\theta}-1)C_4}}{2C_4}, \quad \theta\leq \theta_{NI}\label{eq:s-outer} ,
\end{equation}
with $\bar\theta$ to be determined from the $\theta$ solution of the $\theta$-equation given next.
The initial layer term satisfies
the equation
\begin{equation}
    \odv{S_I}{\tau}=-a_1 S_I - b_1 S_I^2 - c_1 S_I^3, \quad \text {such that}\quad S_I=O(e^{-\frac{\bar t}{C_1}}).\label{eq:s-inner-Neumann}
\end{equation}

We now turn to finding the inner and outer solutions of the temperature equation \eqref{eq:temp_nondim}.  We first observe that the initial layer term associated with solutions of such an equation has the size $O(\frac{\bar t}{B_1})$ about $10^5$ times  larger than that corresponding to the $s$-equation. This together with the presence of the heat operator in the linear part of such  equation allows as to bypass the explicit separation of the initial layer  term $O(e^{-\lambda_1}{C_1})$, in the solution. Furthermore, the fact that the second derivative operator enters as an order-1 term, no boundary layer component is present in the solution.

We  treat the case when Robin boundary conditions for $\theta$ as in \eqref{eq:bc-theta_neumann}, \eqref{eq:bc-s_neumann}, and \eqref{eq:initial-s_theta} are satisfied at $z=-\frac{1}{2}.$ So, the solution of equation \eqref{eq:temp_nondim} is given by 

\begin{equation}
\label{eq:temperature-solution}
\begin{aligned}
\bar{\theta}(\bar z,\bar t)
={}& \bar{\theta}^{+}-\bar h\,\bar t
+\bar h B_1
\left(
\frac{1}{\mathrm{Bi}}
+\bar z-\frac{\bar z^2}{2}
\right)-\bar h B_1
\sum_{n=1}^{\infty}
A_n
\exp\left(
-\frac{\lambda_n^2\bar t}{B_1}
\right)
\cos\bigl[\lambda_n(1-\bar z)\bigr].
\end{aligned}
\end{equation}

\begin{equation} \nonumber
A_n =
\frac{
\displaystyle
\int_0^1
\left(
\frac{1}{\mathrm{Bi}}+\bar z-\frac{\bar z^2}{2}
\right)
\cos\bigl[\lambda_n(1-\bar z)\bigr]
\,d\bar z
}{
\displaystyle
\int_0^1
\cos^2\bigl[\lambda_n(1-\bar z)\bigr]
\,d\bar z
},
\qquad
\lambda_n\tan\lambda_n=\mathrm{Bi}.
\end{equation}
Here, $A_n$ has been computed from the initial condition $\bar{\theta}(\bar{z},0)=\bar{\theta}^+$.
The last  correction term in \eqref{eq:temperature-solution} gives the contribution of the regular perturbation expression $ \partderiv{s}{t}$ in \eqref{eq:temp_nondim}, which by Duhamel's formula, is found to contribute only in a time layer $O(e^{-t/C_1})$.  

Next, we calculate the front in the plane (t, z) that separates the 
nematic from the isotropic region. For this, we set $\bar\theta=\frac{\theta_{NI}}{\theta^*} $ in the previous equation.   That is, \begin{equation}
\label{eq:temperature_first_mode}
\begin{aligned}
\frac{\theta_{NI}}{\theta^*}=1.0357
={}& \bar{\theta}^{+}-\bar h\,\bar t
+\bar h B_1
\left(
\frac{1}{\mathrm{Bi}}
+\bar z-\frac{\bar z^2}{2}
\right)-\bar h B_1 A_1
\exp\left(
-\frac{\lambda_1^2\bar t}{B_1}
\right)
\cos\bigl[\lambda_1(1-\bar z)\bigr].
\end{aligned}
\end{equation}
\begin{figure}[htbp]
    \centering
    \begin{subfigure}[b]{0.4\textwidth}
        \centering
        \includegraphics[width=\textwidth]{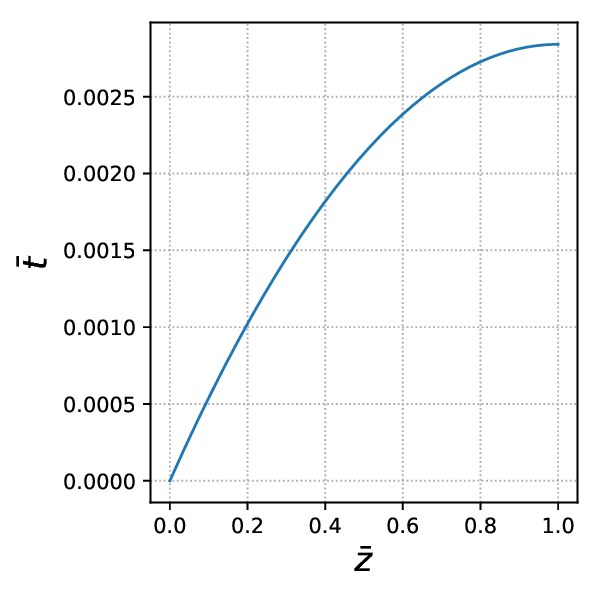}
        \caption{Temperature profiles $\bar{\theta}$}
        \label{fig:temp_profiles}
    \end{subfigure}
    \hfill 
    \begin{subfigure}[b]{0.40\textwidth}
        \centering
        \includegraphics[width=\textwidth]{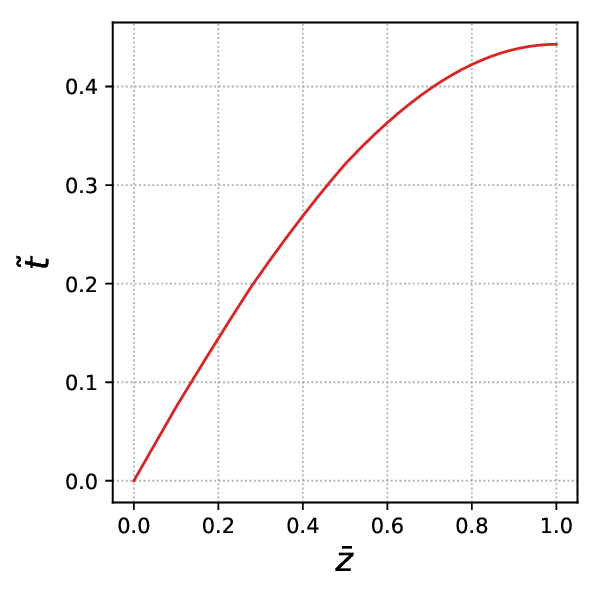}
        \caption{Order parameter profiles $S$}
        \label{fig:order_profiles}
    \end{subfigure}
    \caption{\small{ Plots of the nematic-isotropic front across the dimensionless cylinder  with two different heights and time scales. \text{\bf{(a)}} For the cylinder with $d=0.344 \text{ mm}$, the front delay corresponds to the bath-cooling time-scale $\bar t=\frac{t}{t_\circ}$, $t_\circ=\frac{\text{ K}}{h}\approx 59.86 s. $  \text{\bf{(b)}} 
    For the cylinder with $d=1 \text{ cm}$, the front delay corresponds to the diffusive time-scale $\bar t=\frac{t}{\tau_D}$, $\tau_D\approx 287.5  s. $
    }} 
    \label{fig:combine_profiles}
\end{figure}

\subsection{Liquid crystal diffusion time-scale and asymptotic analysis}
\label{subsec:diffusive_time_scale}

For sample geometries exceeding the critical thickness ($d > d_{\text{crit}} \approx 4.56\text{ mm}$), the external bath cooling time scale $t_0$ loses physical validity as internal heat conduction becomes the rate-limiting mechanism. To analyze the bulk-diffusion regime, we rescale the temporal variable using the characteristic thermal diffusion time scale of the liquid crystal:
\begin{equation}
    \tau := \tau_D = \frac{\rho c_\theta d^2}{\kappa_{||}} = \frac{d^2}{\alpha_T}.
    \label{eq:tau_D_definition}
\end{equation}
Under the temporal transformation $\tilde{t} = t / \tau_D = \bar{t} / B_1$, the non-dimensional system~\eqref{eq:temp_nondim}--\eqref{eq:order_nondim} transforms into:
\begin{align}
    \frac{\partial\bar{\theta}}{\partial\tilde{t}} &= \frac{\partial^2\bar{\theta}}{\partial\bar{z}^2} - \frac{1}{B_1}\left[B_2 s^3 + B_3 s^2 - B_4 s\right] \frac{\partial s}{\partial\tilde{t}}, \label{eq:temp_diff_scaling} \\
    \tilde{C}_1 \frac{\partial s}{\partial\tilde{t}} &= C_2 \frac{\partial^2 s}{\partial\bar{z}^2} - \left[(\bar{\theta}-1)s + C_3 s^2 + C_4 s^3\right], \label{eq:order_diff_scaling}
\end{align}
where $\tilde{C}_1 = C_1 / B_1$. The Robin convective boundary condition at the bottom wall ($\bar{z}=0$) and the insulated Neumann condition at the top wall ($\bar{z}=1$) become[cite: 43]:
\begin{equation}
    \left(\frac{\partial\bar{\theta}}{\partial\bar{z}}\right)_{\bar{z}=0} = Bi \left(\bar{\theta} - \bar{\theta}^+ + \tilde{h}\tilde{t}\right), \qquad \left(\frac{\partial\bar{\theta}}{\partial\bar{z}}\right)_{\bar{z}=1} = 0,
    \label{eq:robin_diff_bc}
\end{equation}
where $\tilde{h} = \bar{h} B_1 = \frac{h \tau_D}{\theta^*}$ represents the dimensionless cooling parameter under the diffusion scaling.

Evaluating equation \eqref{eq:tau_D_definition} with the ITZ physical constants from Table~1 for a sample height of $d = 1.0\text{ cm}$ gives the diffusion time scale $\tau_D \approx 287.5\text{ s}$ ($\approx 4.79\text{ min}$). We calculate the  corresponding dimensionless coefficients \eqref{parameter_groups} to get
\begin{equation}
\begin{aligned}
    B_1 = 4.803, \quad B_2 = 1.85 \times 10^{-3}, \quad   B_3 = -1.70 \times 10^{-3}, \quad  B_4 = 0.116,\quad  Bi = 1.25\\
    \tilde{C}_1 = 2.23 \times 10^{-10},  \quad C_2 = 5.43 \times 10^{-14}, \quad  C_3 = -0.183,  \quad C_4 = 0.200, \quad
    \tilde{h} = 0.0137.
\end{aligned}
\label{eq:diff_coefficients_1cm}
\end{equation}

We note  that while $B_1 \ll 1$ in the previous  geometry, here $B_1 = 4.803 > 1$, indicating the dominant role of the internal heat conduction. Furthermore, the Biot number increases to $Bi = 1.25 \sim O(1)$, indicating that bulk thermal resistance and external convective heat transfer operate on comparable scales.

We now consider the asymptotic approximation that takes into account  that 
 $C_2 \approx 5.43 \times 10^{-14} \ll 1$ and $\tilde{C}_1 \ll 1$ as well. So, the spatial and temporal derivative terms in the order parameter equation~\eqref{eq:order_diff_scaling} are negligible outside a narrow transient initial layer. We carry out the asymptotic study of the equations that result from setting $C_1=0$ and  $C_2 = 0$, and consequently neglecting  the initial and boundary layer contributions in $s$, resulting in the bulk algebraic relation of  the scalar order parameter with respect  to the temperature field as given by \eqref{eq:s-outer}.
Applying the chain rule to the  s-derivative in equation~\eqref{eq:temp_diff_scaling}, $ \frac{\partial s}{\partial\tilde{t}} = \frac{\mathrm{d} S}{\mathrm{d}\bar{\theta}} \frac{\partial\bar{\theta}}{\partial\tilde{t}},$
and factoring out $\frac{\partial\bar{\theta}}{\partial\tilde{t}}$ the system reduces to a single quasilinear heat equation,
\begin{equation}
    C_{\text{eff}}(\bar{\theta}) \frac{\partial\bar{\theta}}{\partial\tilde{t}} = \frac{\partial^2\bar{\theta}}{\partial\bar{z}^2},\qquad   C_{\text{eff}}(\bar{\theta}):= 1 + \frac{1}{B_1} \left[ B_2 S(\bar{\theta})^3 + B_3 S(\bar{\theta})^2 - B_4 S(\bar{\theta}) \right] \frac{\mathrm{d} S}{\mathrm{d}\bar{\theta}},
    \label{eq:eff_heat_pde}
\end{equation}
where $C_{\text{eff}}(\bar{\theta})$ is the effective heat capacity that accounts for the latent heat of the isotropic-nematic transition.

\subsubsection{Front Propagation and Numerical Trajectory}
Let us begin with a heuristic argument that illustrates why a parabolic profile for the front is also to be expected in the diffusion-dominated regime (e.g., $d \ge 1\text{ cm}$). After an initial thermal adjustment, the bulk of the material must cool down at the prescribed  rate of the bath. In the diffusion time scale $\tilde{t}$, this  imposes a uniform dimensionless temperature drop $\frac{\partial\bar{\theta}}{\partial\tilde{t}} \approx -\tilde{h}$ across the sample. Substituting this uniform cooling rate into the dominant heat operator yields a constant spatial curvature, 
$
    \frac{\partial^2\bar{\theta}}{\partial\bar{z}^2} \approx -\tilde{h}.
$
Integrating this balance twice with respect to $\bar{z}$, and applying the insulated top boundary condition ($\left.\frac{\partial\bar{\theta}}{\partial\bar{z}}\right|_{\bar{z}=1} = 0$), yields an exact quadratic spatial temperature distribution,
$
    \bar{\theta}(\bar{z}, \tilde{t}) \approx \bar{\theta}(0, \tilde{t}) + \tilde{h} \left( \bar{z} - \frac{\bar{z}^2}{2} \right).
$
Evaluating the former at the phase transition threshold $\bar{\theta}(\bar{z}, \tilde{t}) = \bar{\theta}_{NI}$ renders the trajectory of the moving phase boundary. By converting the dimensionless cooling parameter $\tilde{h}$ back to physical parameters, we find that the physical time delay $\Delta t(\bar{z}) = t(\bar{z}) - t(0)$ required for the front to travel from the cooled base ($\bar{z}=0$) to coordinate $\bar{z}$ must follow an exact concave-downward parabola:
\begin{equation}
    \Delta t(\bar{z}) \approx B_1 \tau_D \left( \bar{z} - \frac{\bar{z}^2}{2} \right).
    \label{eq:asymptotic_parabola_delay}
\end{equation}
Physically, because heat can only leave through the bottom boundary ($\bar{z}=0$), the front advances rapidly at early times near the base where thermal gradients are steep, but  slows down as it approaches the insulated top wall where the local thermal gradient gets close to zero, that is, $\left.\frac{\partial(\Delta t)}{\partial\bar{z}}\right|_{\bar{z}=1} \to 0$.
\subsubsection{Numerical Scheme for the Effective Quasilinear Equation}
To simulate the full non-linear front propagation without relying on quasi-steady asymptotic approximations, equation~\eqref{eq:eff_heat_pde} was integrated numerically using the Method of Lines (MoL) \cite{schiesser1991numerical}:

\begin{enumerate}
    \item \textbf{Spatial Discretization:} The domain $\bar{z} \in [0, 1]$ is discretized into $N = 150$ uniform mesh nodes with grid spacing $\Delta \bar{z} = \frac{1}{N-1}$. The spatial Laplacian is discretized using a second-order centered finite difference operator:
    \begin{equation}
        \left. \frac{\partial^2 \bar{\theta}}{\partial \bar{z}^2} \right|_{\bar{z}_i} \approx \frac{\bar{\theta}_{i+1} - 2\bar{\theta}_i + \bar{\theta}_{i-1}}{(\Delta \bar{z})^2}.
    \end{equation}

    \item \textbf{Boundary Condition Implementation:} To preserve global second-order spatial accuracy, boundary conditions are enforced using virtual ghost nodes $\bar{\theta}_{-1}$ and $\bar{\theta}_{N}$:
    \begin{itemize}
        \item \emph{Cooled Bottom Wall ($\bar{z} = 0$):} The Robin convective condition $\left.\frac{\partial \bar{\theta}}{\partial \bar{z}}\right|_{\bar{z}=0} = Bi \left(\bar{\theta}_0 - \bar{\theta}_b(\tilde{t})\right)$ specifies $\bar{\theta}_{-1} = \bar{\theta}_1 - 2 \Delta \bar{z} Bi \left(\bar{\theta}_0 - \bar{\theta}_b(\tilde{t})\right)$.
        \item \emph{Insulated Top Wall ($\bar{z} = 1$):} The Neumann condition $\left.\frac{\partial \bar{\theta}}{\partial \bar{z}}\right|_{\bar{z}=1} = 0$ sets $\bar{\theta}_{N} = \bar{\theta}_{N-2}$.
    \end{itemize}

    \item \textbf{Temporal Integration:} The resulting system of $N$ coupled non-linear ODEs,
    \begin{equation}
        \frac{\mathrm{d} \bar{\theta}_i}{\mathrm{d} \tilde{t}} = \frac{1}{C_{\mathrm{eff}}(\bar{\theta}_i)} \left[ \frac{\bar{\theta}_{i+1} - 2\bar{\theta}_i + \bar{\theta}_{i-1}}{(\Delta \bar{z})^2} \right],
    \end{equation}
    is integrated using an implicit stiff Radau IIA solver over dimensionless time $\tilde{t} \in [0, 15.0]$ with strict error tolerances ($\mathrm{rtol} = 10^{-8}$, $\mathrm{atol} = 10^{-10}$) \cite{hairer1996solving}.

\item \textbf{Temporal Interpolation for Isotherm Tracking } \cite{crank1984free, ShampineReichelt1997, OsherSethian1988}:
To determine the propagation of the transition-temperature isotherm \(\bar{\theta}=\bar{\theta}_{NI}\), we calculate the first crossing time at each spatial node \(\bar{z}_i\).
Suppose that the temperature crosses \(\bar{\theta}_{NI}\) between two consecutive computed time levels \(\tilde{t}_n\) and \(\tilde{t}_{n+1}\), so that
\[
\bar{\theta}_i^n > \bar{\theta}_{NI}
\geq \bar{\theta}_i^{n+1},
\qquad
\bar{\theta}_i^n
:= \bar{\theta}(\bar{z}_i,\tilde{t}_n).
\]
The crossing time is then approximated by linear interpolation in time:
\begin{equation}
\tilde{t}_{\mathrm{front}}(\bar{z}_i)
\tilde{t}n
+
\frac{\bar{\theta}{NI}-\bar{\theta}_i^n}
{\bar{\theta}_i^{n+1}-\bar{\theta}i^n}
(\tilde{t}{n+1}-\tilde{t}_n).
\label{eq:isotherm_crossing_time}
\end{equation}
This procedure reduces the temporal discretization noise associated with assigning the crossing to the nearest computed time level.
The corresponding physical propagation delay, measured relative to the bottom boundary, is
\begin{equation}
\Delta t(\bar{z}_i)
\tau_D
\left[
\tilde{t}_{\mathrm{front}}(\bar{z}_i)
\tilde{t}_{\mathrm{front}}(0)
\right].
\label{eq:isotherm_propagation_delay}
\end{equation}
\end{enumerate}
For $d = 1.0~\mathrm{cm}$, the full PDE numerical solution predicts a total traversal duration across the cylinder height of $\Delta t_{\mathrm{max}} \approx 143.8~\mathrm{s}$ ($\approx 2.4~\mathrm{min}$), confirming the dominant role of internal thermal conduction in thick sample processing.

\begin{figure}[htbp]
    \centering
    \begin{subfigure}[b]{0.4\textwidth}
        \centering
        \includegraphics[width=\textwidth]{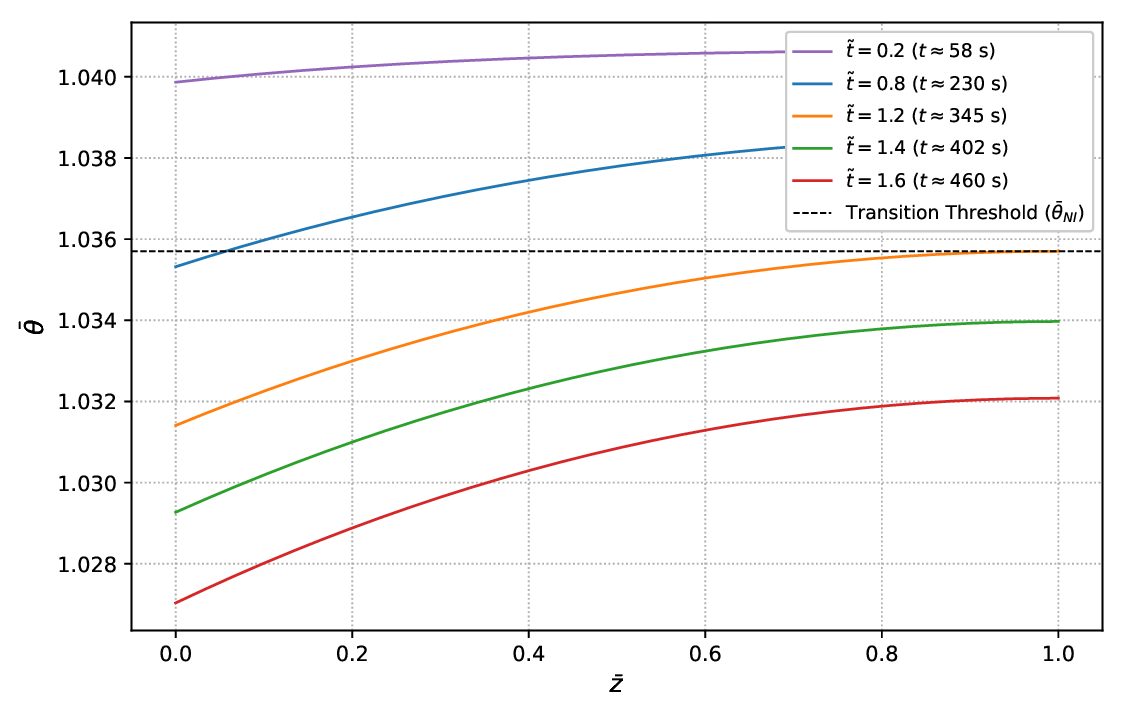}
        \caption{Temperature profiles $\bar{\theta}$}
        \label{fig:temp_profiles}
    \end{subfigure}
    \hfill 
    \begin{subfigure}[b]{0.40\textwidth}
        \centering
        \includegraphics[width=\textwidth]{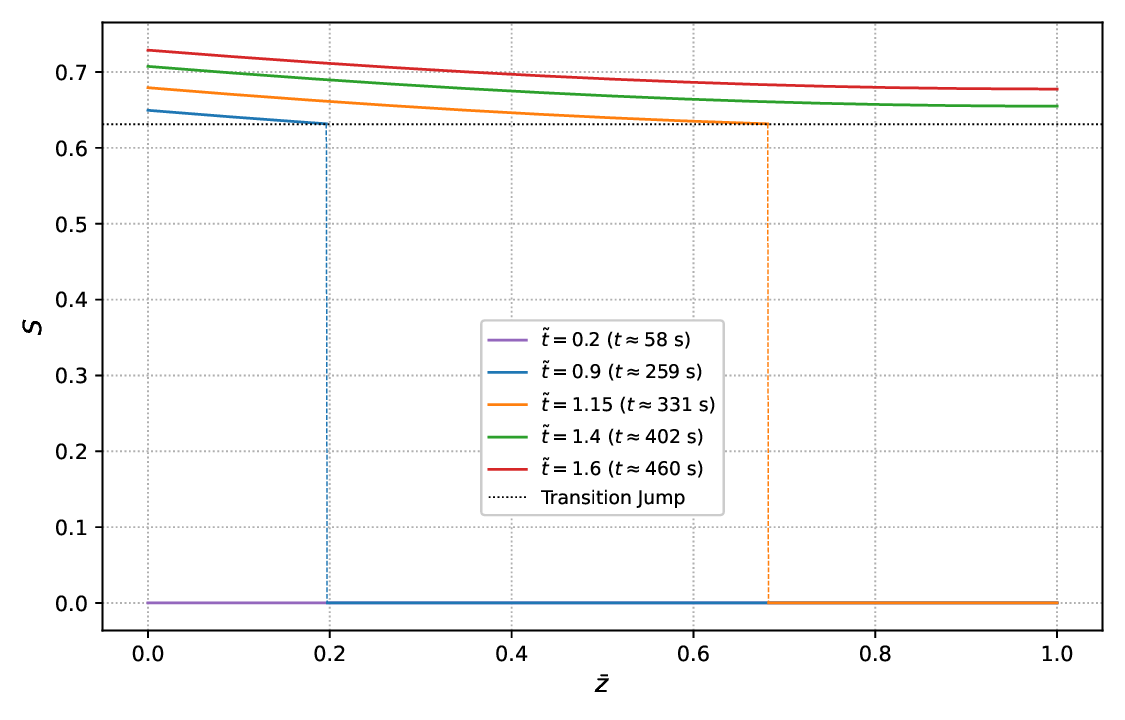}
        \caption{Order parameter profiles $S$}
        \label{fig:order_profiles}
    \end{subfigure}
    \caption{\small{Spatial profiles of the temperature (a) and scalar order parameter (b) across a $1.0\text{ cm}$ sample at successive dimensionleYourss times. The profile at $\tilde{t}=0.25$ show the sample in the isotropic phase prior to nucleation. As the bulk cools, the discontinuous isotropic-nematic phase front  sweeps across the domain ($\tilde{t}=0.25, \,\tilde{t}=0.25$).  The final profiles ($\tilde{t}=1.45$ and  \,$\tilde{t}=1.65$) show the nematic order steadily increasing as the front  reaches the insulated top boundary.  } }
    \label{fig:combined_profiles}
\end{figure}

\section {Conclusions}

We have investigated the cooling-induced isotropic–nematic phase transition in thermotropic liquid crystals, motivated by the processing of liquid-crystalline pharmaceutical materials and, in particular, itraconazole. Within a non-isothermal extension of Ericksen's theory of liquid crystals with variable degree of orientation, we have developed a thermodynamically consistent description of the coupling between orientational relaxation and heat conduction. The resulting nonlinear parabolic system accounts for the evolution of the scalar order parameter and temperature, including the effects of anisotropic thermal conductivity and heat exchange with an externally controlled cooling bath.

We highlight two aspects of the mathematical formulation. First, along with the classical polynomial Landau–de Gennes bulk potential, we have considered the singular Ball–Majumdar potential, which incorporates the physical restrictions on the eigenvalues of the orientational order tensor. To overcome the implicit character of this potential, we have constructed an explicit approximation that retains its endpoint singularities and the coexistence conditions associated with the isotropic–nematic transition. We have established convergence of this approximation and of the corresponding coexistence parameters. Second, we have investigated the well-posedness of the coupled non-isothermal problem, establishing local existence, uniqueness, and regularity of solutions, together with separation from the physical bounds of the order parameter and a continuation criterion.

The nondimensional and asymptotic analysis reveals the importance of the widely separated time scales governing orientational relaxation, thermal diffusion, and externally imposed cooling. In particular, comparison of the thermal diffusion time with the cooling time identifies a characteristic sample thickness separating two regimes of thermal evolution. For sufficiently thin samples, the temperature remains nearly uniform across the material, whereas for thicker samples the finite rate of thermal diffusion produces significant spatial delays. The reduced thermal equations yield approximate descriptions of the propagation of the isotropic–nematic transition-temperature isotherm, with a concave quadratic dependence of its time of evolution  on the vertical coordinate. These results illustrate how changes in sample dimensions may substantially modify the thermal history during cooling. The transition-temperature isotherm should, however, be distinguished from the actual isotropic–nematic interface, whose evolution also depends on the kinetics of orientational ordering and possible metastability.

Several extensions of the present work are of interest. An extension of the model to include the subsequent nematic–smectic A transition is currently under way. This requires introducing additional order parameters to describe the development of smectic layering and their coupling to nematic orientation and temperature. Such an extension is particularly relevant to itraconazole, for which the sequence of isotropic, nematic, and smectic phases plays an important role in determining the material state attained during cooling. Other directions include a more detailed investigation of the dynamics of the phase interface, the influence of temperature-dependent material coefficients and nucleation mechanisms. These developments would provide a basis for studying the competition between phase ordering, thermal transport, and kinetic arrest in the preparation of liquid-crystalline pharmaceutical materials.


\bibliographystyle{siam}
\bibliography{LC}

\end{document}